\documentclass{elsarticle}

\usepackage{graphicx}
\usepackage{sidecap}
\usepackage{wrapfig}
\usepackage{amsmath,amsthm,amsfonts,amssymb,amscd,mathrsfs,eucal}
\usepackage{fancyhdr}

\def\r{\rho}
\def\d{\delta}
\def\e{\epsilon}
\def\g{\gamma}
\def\a{\alpha}

\def\lmax{\lambda^{max}}
\def\lmin{\lambda^{min}}

\newtheorem{theorem}{Theorem}
\newtheorem{corollary}[theorem]{Corollary}

\newtheorem{definition}[theorem]{Definition}

\journal{Linear Algebra and its Applications}

\begin{document}

\begin{frontmatter}

\title{Bounds of restricted isometry constants in extreme asymptotics: formulae for Gaussian matrices}

\author[bb]{Bubacarr˜Bah\corref{cor1}}
\ead{b.bah@sms.ed.ac.uk}
\author[jt]{Jared˜Tanner\fnref{fn1}}
\ead{tanner@maths.ox.ac.uk}

\cortext[cor1]{Corresponding author}
\fntext[fn1]{This author's work was supported in part by the Leverhulme Trust.}

\address[bb]{Maxwell Institute and School of Mathematics, University of Edinburgh, Edinburgh, UK}
\address[jt]{Mathematics Institute and Exeter College, University of Oxford, Oxford, UK}

%\baselineskip=20pt plus1pt
%\maketitle
%\tableofcontents
%\pagestyle{headings}
%\markright{Bounds of restricted isometry constants in extreme asymptotics}

\begin{abstract}
{\em Restricted Isometry Constants (RICs) provide a measure of how far from an isometry a matrix can be when acting on sparse vectors.  This, and related quantities, provide a mechanism by which standard eigen-analysis can be applied to topics relying on sparsity.  RIC bounds have been presented for a variety of random matrices and matrix dimension and sparsity ranges.  We provide explicitly formulae for RIC bounds, of $n \times N$ Gaussian matrices with sparsity $k$, in three settings: a) $n/N$ fixed and $k/n$ approaching zero, b) $k/n$ fixed and $n/N$ approaching zero, and c) $n/N$ approaching zero with $k/n$ decaying inverse logrithmically in $N/n$; in these three settings the RICs a) decay to zero, b) become unbounded (or approach inherent bounds), and c) approach a non-zero constant.  Implications of these results for RIC based analysis of compressed sensing algorithms are presented.}
\end{abstract}


%\subclass{Primary, 15B52, 60F10, 94A20; Secondary, 94A12}

\end{frontmatter}

\section{Introduction}\label{sec:Intro}
Many questions in signal processing\cite{baraniuk2011more,haupt2007compressive}, statistics \cite{babacan2010bayesian,davenport2006detection,lounici2009taking}, computer vision \cite{cevher2008sparse,romberg2008imaging,stankovic2008compressive,wright2010sparse}, and machine learning \cite{calderbank2009compressed,cevher2008learning,mahoor2011facial} are employing a parsimonious notion of eigen-analysis to better capture inherent simplicity in the data. Slight variants of the same quantity are defined in these disciplines, referred to as: sparse principal components, sparse eigenvalues, and restricted isometry constants (RICs).  In this article we adopt the notation and terminology of RICs, defined as a measure of the greatest relative change that a matrix can induce in the $\ell^2$ norm of sparse vectors.  Let $\chi^N(k)$ denote all vectors of length $N$ which have at most $k$ nonzeros; then the lower and upper RICs of the $n\times N$ matrix $A$ are defined as
\begin{eqnarray}
L(k,n,N;A) & := & 1-\min_{x\in\chi^N(k)} \frac{\|Ax\|_2^2 }{ \|x\|_2^2 }
\quad \mbox{and}\label{eq:L} \\
U(k,n,N;A) & := & \max_{x\in\chi^N(k)} \frac{\|Ax\|_2^2 }{ \|x\|_2^2 } - 1
\quad\mbox{respectively.}\label{eq:U}
\end{eqnarray}
RICs were introduced by Cand\`es and Tao in 2004 \cite{candes2005decoding} as a method of analysis for sparse approximation and compressed sensing (CS), and have received widespread used in those communities.  For example, let $y=Ax_0+e$ for some $x_0\in\chi^N(k)$, then, provided the RICs of $A$ are sufficiently small, there are computationally tractable algorithms which from $A$ and $y$ (and possibly $k$ and $\|e\|$) are guaranteed to return a vector ${\hat x}$ satisfying a bound of the form $\|x_0-{\hat x}\|_2\le Const. \|e\|_2$; for examples of such theorems see \cite{blanchard2011phase,blumensath2009iterative,candes2008restricted,foucart2009sparsest,dai2009subspace,needell2009cosamp}.  The efficacy of theorems of this form depends highly on knowledge of the RICs of $A$.

Numerous algorithms exist for estimating or bounding the RICs of a general matrix; however, theory for the current state of the art \cite{daspremont2011testing,juditsky2011verifiable} is limited to $k\sim \sqrt{n}$, whereas many applications require information for comparatively larger values of $k$. The only method for calculating the RICs of a general matrix $A$ for larger values of $k$, requires calculating the extreme singular values of all $\binom{N}{k}$ submatrices of $A$, resulting from all independent selections of $k$ columns from $A$.  This combinatorial approach is intractable for all but very small dimensions.  For this reason, much of the research on RICs has been devoted to deriving their bounds.  Matrices with entries drawn from the Gaussian distribution ${\cal N}\left(0,1/n\right)$ have the smallest known bound for large matrices and $k\gg 1$ \cite{bah2010improved}.  For bounds on the RICs of matrix ensembles other than Gaussian see \cite{bajwa2007toeplitz,rauhut2010compressive,fornasier2010theoretical}.

Let
\[
\rho_n:= \frac{k}{n}\quad\mbox{ and }\quad \delta_n:=\frac{n}{N}.
\]
RIC bounds for Gaussian matrices have been derived focusing on the limits $\rho_n\rightarrow\rho\in(0,1)$ and $\delta_n\rightarrow\delta\in(0,1)$, \cite{bah2010improved,blanchard2011compressed,candes2005decoding}, see Theorem \ref{thm:LUbounds}. Unfortunately, these bounds are given in terms of implicitly defined functions, Definition \ref{def:LUbounds}, obscuring their dependence on $\rho$ and $\delta$.

\begin{theorem}[Gaussian RIC Bounds \cite{blanchard2011compressed}]\label{thm:LUbounds}
Let $\mathcal{L}\left(\d,\r\right)$ and $\mathcal{U}\left(\d,\r\right)$ be defined as in Definition \ref{def:LUbounds} and fix $\e>0$. In the limit where $n/N\rightarrow\d\in (0,1)$ and $k/n\rightarrow\r\in (0,1)$ as $n\rightarrow\infty$, sample each $n\times N$ matrix $A$ from the Gaussian ensemble (entries drawn independent and identically distributed from the Gaussian Normal ${\cal N}\left(0,1/n\right)$) then
\begin{align*}
& \hbox{Prob}\left(L(k,n,N;A)<\mathcal{L}\left(\d,\r\right)+\e\right)\rightarrow 1\quad \hbox{and}\quad \\
& \hbox{Prob}\left(U(k,n,N;A)<\mathcal{U}\left(\d,\r\right)+\e\right)\rightarrow 1
\end{align*}
exponentially in $n$.
\end{theorem}

In this manuscript we present simple expressions which bound the RICs of Gaussian matrices in three asymptotic settings: (a) $\d\in (0,1)$ and $\r\ll 1$ where the RICs converge to zero as $\r$ approaches zero, (b) $\r\in (0,1)$ and $\d\ll 1$ where the upper RIC become unbounded and the lower RIC converges to its bound of one as $\d$ approaches zero, and (c) along the path $\r_{\gamma}(\d)=1/\left(\gamma\log(\d^{-1})\right)$ for $\d\ll 1$ where the RICs approach a nonzero constant as $\d$ approaches zero. In all cases, except for the bound of the lower RIC in case b) we see the introduction of a new logarithmic term coming from the combinatorial term which is a result of the union bound we use in the derivations (see proof of the main results). Furthermore, we have a $\d$ dependence in the factor $\d^2\r^3$ in all the bounds.

The bounds presented here build on the results  in \cite{blanchard2011compressed} and are specific to Gaussian matrices, carefully balancing combinatorial quantities with the tail behaviour of the largest and smallest singular values of Gaussian matrices. The specificity of these bounds to Gaussian matrices gives great accuracy than what subgaussian tail bounds provide \cite{baraniuk2008simple}. A similar analysis could be conducted for the subgaussian case by considering the bounds in \cite{candes2005decoding} stated for the Gaussian case, but which are equally valid for the subgaussian case.  For brevity we do not consider the subgaussian case here.

There has been substantial work on RICs of partial  Fourier matrices, see \citep{rauhut2010compressive} and references therein. However, the exact power of the logarithmic factor (in  $(\gamma\log(1/\delta_n))^{-1}$) is not yet determined. Hence analysis of the kind of this work are not possible for such ensembles.

Each of Theorems \ref{thm:ULbctr} -- \ref{thm:ULbctg} state that the probability under consideration {\em converge exponentially} to $1$ in $k$ or $n$ which we use as a shorthand for saying one minus the probability considered being bounded by a function decaying exponentially to zero in the variable stated; the explicit bound is given in the proof of the theorem.

Theorem \ref{thm:LUbounds} states that, for $k$, $n$, and $N$ large, it is unlikely that the RICs exceed the constants $\mathcal{L}\left(\d,\r\right)$ and $\mathcal{U}\left(\d,\r\right)$ by more than any $\epsilon$.  In the limit where $\delta_n\rightarrow\delta\in(0,1)$ and $\rho_n\rightarrow \r\ll 1$, the matrix RICs converge to zero, causing the resulting bounds to become vacuous.  Theorem \ref{thm:ULbctr} states the dominant terms in the bounds, and that the true RICs are unlikely to exceed these bounds by a multiplicative factor $(1+\epsilon)$ for any $\epsilon>0$.  The dominant terms can be contrasted with $2\sqrt{\r} + \r$ which is the deviation from one of the expected value of the smallest and largest eigenvalues of a Wishart matrix \cite{Geman1980limit,silverstein1985smallest}.  An implication of Theorem \ref{thm:ULbctr} for the compressed sensing algorithm Orthogonal Matching Pursuit is given in Corollary \ref{omp_corl}.

\begin{theorem}[Gaussian RIC Bounds: $\r\ll 1$]\label{thm:ULbctr}
Let $\widetilde{\mathcal{U}}^{\r}(\delta,\rho)$ and $\widetilde{\mathcal{L}}^{\r}(\delta,\rho)$ be defined as
\begin{align}
\label{eq:Ubctr}
\widetilde{\mathcal{U}}^{\r}(\delta,\rho) & = \sqrt{2\rho\log\left(\frac{1}{\delta^2\rho^3}\right) + c\rho},\\
\label{eq:Lbctr}
\widetilde{\mathcal{L}}^{\r}(\delta,\rho) & = \sqrt{2\rho\log\left(\frac{1}{\delta^2\rho^3}\right) + c\rho}.
\end{align}
Fix $\e>0$ and $c>6$. For each $\d\in(0,1)$ there exists a $\r_0>0$ such that in the limit where $n/N\rightarrow\d$, $k/n\rightarrow\r\in (0,\rho_0)$, and $(\log n)/k\rightarrow 0$ as $k\rightarrow\infty$, sample each $n\times N$ matrix $A$ from the Gaussian ensemble, ${\cal N}\left(0,1/n\right)$, then
\begin{align*}
&\hbox{Prob}\left(L(k,n,N;A)<(1+\e)\widetilde{\mathcal{L}}^{\r}\left(\d,\r\right)\right)\rightarrow 1\quad \hbox{and}\quad \\
&\hbox{Prob}\left(U(k,n,N;A)<(1+\e)\widetilde{\mathcal{U}}^{\r}\left(\d,\r\right)\right)\rightarrow 1
\end{align*}
exponentially in $k$.
\end{theorem}

Theorem \ref{thm:ULbctd} considers a limiting case where the upper RIC diverges and the lower RIC converges to its bound of one. The upper RIC is shown to grow in this setting with a dominant term proportional to $\log(1/\delta)$ with precise proportionality constants as well as the secondary growth factor $\log\log(1/\delta)$, again with constants of proportionality.  The lower RIC is shown to differ from the unit bound by a polynomial term in $\delta$, as opposed to the more typical logarithmic relations.  The rapid decay to zero of the $\delta$ polynomial term in \eqref{eq:Lbctd} indicates that the lower RIC rapidly approaches one as $\delta$ decreases for $\rho$ fixed; this is reflected in the dominant effect of the lower RIC when used to prove convergence guarantees for sparse approximation algorithms \cite{blanchard2011phase}.

%With the lower RIC converging to one, its bound is shown to be no more than an arbitrarily small multiplicative constant, whereas the upper RIC is bounded by an additive constant.

\begin{theorem}[Gaussian RIC Bounds: $\d\ll 1$]\label{thm:ULbctd}
Let $\widetilde{\mathcal{U}}^{\d}(\delta,\rho)$ and $\widetilde{\mathcal{L}}^{\d}(\delta,\rho)$ be defined as
\begin{eqnarray}
\label{eq:Ubctd}
\widetilde{\mathcal{U}}^{\d}(\delta,\rho) & = & \rho\log\left(\frac{1}{\delta^2\rho^3}\right) + (1+\rho)\log\left[c\log\left(\frac{1}{\delta^2\rho^3}\right)\right] + 3\rho, \qquad \\
\label{eq:Lbctd}
\widetilde{\mathcal{L}}^{\d}(\delta,\rho) & = & 1 - \exp\left(-\frac{3\rho+c}{1-\rho}\right) \cdot \left(\delta^2\rho^3\right)^{\frac{\rho}{1-\rho}}.
\end{eqnarray}
Fix $\e>0$ and $c>1$.  For each $\r\in(0,1)$ there exists a $\d_0>0$ such that in the limit where $k/n\rightarrow\r$, $n/N\rightarrow\d\in (0,\d_0)$ as $n\rightarrow\infty$, sample each $n\times N$ matrix $A$ from the Gaussian ensemble, ${\cal N}\left(0,1/n\right)$, then
\begin{align*}
&\hbox{Prob}\left(L(k,n,N;A)<(1+\e)\widetilde{\mathcal{L}}^{\d}\left(\d,\r\right)\right)\rightarrow 1\quad \hbox{and}\quad \\
&\hbox{Prob}\left(U(k,n,N;A)<\widetilde{\mathcal{U}}^{\d}\left(\d,\r\right)+\e\right)\rightarrow 1
\end{align*}
exponentially in $n$.
\end{theorem}

Theorem \ref{thm:ULbctg} considers the path in which both $\rho_n$ and $\delta_n$ converge to zero, but in such a way that the RICs approach nonzero constants.  This path is of particular interest in applications where RICs are required to remain bounded, but where the most extreme advantages of the method are achieved for one of the quantities approaching zero.  For example, compressed sensing achieves increased gains in undersampling as $\delta_n$ decreases to zero; however, all compresses sensing algorithmic guarantees involving RICs require the RICs to remain bounded.  The limit considered in Theorem \ref{thm:ULbctg} provides explicit formula for these algorithms in the case where the undersampling is greatest, see Corollary \ref{alg_corl}.

\begin{theorem}[Gaussian RIC Bounds: $\rho_n\rightarrow (\gamma\log(1/\delta_n))^{-1}$ and $\d\ll 1$]\label{thm:ULbctg}
Let $\rho_{\gamma}(\delta) = \frac{1}{\gamma\log\left(\delta^{-1}\right)}$ and let $\widetilde{\mathcal{U}}^{\gamma}\left(\delta,\rho_{\gamma}(\delta)\right)$ and $\widetilde{\mathcal{L}}^{\gamma}\left(\delta,\rho_{\gamma}(\delta)\right)$ be
defined as
\begin{multline}
\label{eq:Ubctg}
\widetilde{\mathcal{U}}^{\gamma}\left(\delta,\rho_{\gamma}(\delta)\right) = \sqrt{2\rho_{\gamma}(\delta)\log\left(\frac{1}{\delta^2\rho_{\gamma}^3(\delta)}\right) + 6\rho_{\gamma}(\delta)} \\ + c_u\left[2\rho_{\gamma}(\delta)\log\left(\frac{1}{\delta^2\rho_{\gamma}^3(\delta)}\right) + 6\rho_{\gamma}(\delta)\right]
\end{multline}
\begin{multline}
\label{eq:Lbctg}
\widetilde{\mathcal{L}}^{\gamma}\left(\delta,\rho_{\gamma}(\delta)\right) = \sqrt{2\rho_{\gamma}(\delta)\log\left(\frac{1}{\delta^2\rho_{\gamma}^3(\delta)}\right) + 6\rho_{\gamma}(\delta)} \\ - c_l\left[\rho_{\gamma}(\delta)\log\left(\frac{1}{\delta^2\rho_{\gamma}^3(\delta)}\right)  + 6\rho_{\gamma}(\delta)\right].
\end{multline}
Fix $\gamma>\gamma_0$ (which $\gamma_0 \geq 4$), $\e>0$, $c_u>1/3$ and $c_l<1/3$.  There exists a $\d_0>0$ such that in the limit where $k/n\rightarrow\r_{\gamma}(\delta_0)$, $n/N\rightarrow\d\in(0,\d_0)$ as $n\rightarrow\infty$, sample each $n\times N$ matrix $A$ from the Gaussian ensemble, ${\cal N}\left(0,1/n\right)$, then
\begin{align*}
&\hbox{Prob}\left(L(k,n,N;A)<\widetilde{\mathcal{L}}^{\gamma}\left(\d,\rho_{\gamma}(\delta)\right)+\e\right)\rightarrow 1\quad \hbox{and}\quad \\
&\hbox{Prob}\left(U(k,n,N;A)<\widetilde{\mathcal{U}}^{\gamma}\left(\d,\rho_{\gamma}(\delta)\right)+\e\right)\rightarrow 1
\end{align*}
exponentially in $n$.
\end{theorem}

Theorem \ref{thm:ULbctg} considers the path $\r_{\gamma}(\d)$ for $\d\ll 1$; passing to the limit of $\d\rightarrow 0$, the functions
$\widetilde{\mathcal{U}}^{\gamma}\left(\delta,\rho_{\gamma}(\delta)\right)$ and $\widetilde{\mathcal{L}}^{\gamma}\left(\delta,\rho_{\gamma}(\delta)\right)$
defined as \eqref{eq:Ubctg} and \eqref{eq:Lbctg} converge to simple functions of $\gamma$.

\begin{corollary}[Gaussian RIC Bounds: $\rho_n\rightarrow (\gamma\log(1/\delta_n))^{-1}$ as $\d\rightarrow 0$]\label{cor_ULrgd}
Let $\widetilde{\mathcal{U}}^{\gamma}\left(\delta,\rho_{\gamma}(\delta)\right)$ and $\widetilde{\mathcal{L}}^{\gamma}\left(\delta,\rho_{\gamma}(\delta)\right)$ be
defined as \eqref{eq:Ubctg} and \eqref{eq:Lbctg} respectively with $\r_{\gamma}(\delta)=\frac{1}{\gamma\log(\d^{-1})}$.
\begin{eqnarray}
\label{cor_Urgd}
\lim_{\delta\rightarrow 0} \widetilde{\mathcal{U}}^{\gamma}\left(\delta,\rho_{\gamma}(\delta)\right) & = & \frac{2}{\sqrt{\gamma}} + \frac{4}{\gamma}c_u \\
\label{cor_Lrgd}
\lim_{\delta\rightarrow 0} \widetilde{\mathcal{L}}^{\gamma}\left(\delta,\rho_{\gamma}(\delta)\right) & = & \frac{2}{\sqrt{\gamma}} - \frac{4}{\gamma}c_l.
\end{eqnarray}
\end{corollary}

The accuracy of Theorems \ref{thm:ULbctr} - \ref{thm:ULbctg} and Corollary \ref{cor_ULrgd} are discussed in Section \ref{sec:rip_asymp} and proven in Section \ref{sec:prip_asymp}.

\subsection{Compressed sensing sampling theorems}\label{sec:mresults}

Compressed sensing is a technique by which simplicity in data can be exploited to reduce the amount of measurements needed to acquire the data.  For example, let there be a vector $x_0\in\chi^N(k)$ which satisfies $y=Ax_0+e$; the matrix $A$ can be viewed as measuring $x_0$ through inner products between its rows and $x_0$, and $e$ captures the model misfit such as measurement error or the true measured vector not being exactly $k$ sparse.  If we let $A$ be of size $n\times N$ with $n<N$, then fewer than $N$ inner products have been performed, and naively it seems impossible to recover $x_0$.

The theory of compressed sensing has developed conditions in which $x_0$, or an approximation thereof, can be recovered.  Most remarkably, for any fixed ratio $n/N$, the recovery guarantees achieve the optimal order of the number of measurements being proportional to the information content in $x_0$ ($n$ proportional to $k$).  In fact, for most compressed sensing algorithms it is possible to derive constants of proportionality, $\rho^{alg}(\delta)$, such that if $A$ has entries ${\cal N}\left(0,1/n\right)$, then in the limit of $n\rightarrow\infty$ with $n/N\rightarrow\delta\in (0,1)$ and $k/n<(1-\epsilon)\rho^{alg}(\delta)$ it can be guaranteed that the output of a compressed sensing algorithm, ${\hat x}$, will satisfy $\|x_0-{\hat x}\|_2\le Const.\|e\|_2$.  The best current known values of $\rho^{alg}(\delta)$ have been calculated in \cite{blanchard2011phase} for Iterative Hard Thresholding (IHT) \cite{blumensath2009iterative}, Subspace Pursuit (SP) \cite{dai2009subspace}, and Compressed Sampling Matching Pursuit (CoSaMP) \cite{needell2009cosamp}. It can be expected that further analysis of these algorithms will result in higher phase transitions, $\rho^{alg}(\delta)$.

Compressed sensing is most remarkable in that the recovery algorithms remain effective for $k/n$ decaying slowly as the number of measurements becomes vanishingly small compared to the signal length, $n/N\rightarrow 0$.  In fact, it is known that $\rho^{alg}(\delta)$ becomes proportional to $1/\left(\log(\delta^{-1})\right)$ as $\delta\rightarrow 0$.  This constant of proportionality can be deduced from Theorem \ref{thm:ULbctg}; the resulting sampling theorems for representative compressed sensing algorithms are stated in Corollary \ref{alg_corl} for $c_u=c_l=1/3$.

\begin{corollary}\label{alg_corl}
Given a sensing matrix, $A$, of size $n\times N$ whose entries are drawn i.i.d. from $\mathcal{N}\left(0,1/n\right)$, in the limit as $n/N \rightarrow 0$ a sufficient condition for recovery for Compressed Sensing algorithms with $k$ steps is $n \ge \gamma k \log\left(N/n\right)$ measurements with $\gamma = 37$ for $l_1$-minimization \cite{candes2008restricted}, $\gamma = 96$ for Iterative Hard Thresholding (IHT) \cite{blumensath2009iterative}, $\gamma = 279$ for Subspace Pursuit (SP) \cite{dai2009subspace}, and $\gamma = 424$ for Compressed Sampling Matching Pursuit (CoSaMP) \cite{needell2009cosamp};  while $\gamma = 332$ for Orthogonal Matching Pursuit (OMP) with $31k$ steps \cite{zhang2011sparse}.
\end{corollary}

Not all compressed sensing algorithms achieve the optimal order of $k$ being proportional to $n$ with $k$ steps. That is converging, to the exact solution for the noiseless case or to the desired approximation error when the measurements have noise, after $k$ steps with the number of measurements $n$ being proportional to $k$, i.e. $n = \mathcal{O}(k\log(N/k))$. One such algorithm is Orthogonal Matching Pursuit (OMP), which has recently been analyzed using RICs, see \cite{zhang2011sparse,mo2011remarks} and references therein.  An analytic asymptotic sampling theorem for OMP with $k$ steps can be deduced from Theorem \ref{thm:ULbctr}, see Corollary \ref{omp_corl}.

\begin{corollary}\label{omp_corl}
Given a sensing matrix, $A$, of size $n\times N$ whose entries are drawn i.i.d. from $\mathcal{N}\left(0,1/n\right)$, in the limit as $n/N \rightarrow\d\in(0,1)$ a sufficient condition for recovery for Orthogonal Matching Pursuit (OMP) with $k$ steps is
\[n> 2k(k-1)[3 + 2\log N + \log n - 3\log k].\]
\end{corollary}

%However, OMP has been proven to give exact recovery of $k$ sparse vectors when $31k$ steps are taken, see \cite{zhang2011sparse} with a condition equivalent to saying the lower RIC to be less than $1/3$, and that in this setting OMP would have a recovery similar to those stated in Corollary \ref{alg_corl}, but with $\gamma = 332$.

\section{Accuracy of main results}\label{sec:rip_asymp}
This section discusses the accuracy of Theorems \ref{thm:ULbctr} - \ref{thm:ULbctg} and Corollary \ref{cor_ULrgd}, comparing the expressions with the bounds in Theorem \ref{thm:LUbounds}, which are defined \cite{blanchard2011compressed} implicitly in Definition \ref{def:LUbounds}.

\begin{definition}\label{def:LUbounds}
Define $\mathcal{L}(\d,\r)$ and $\mathcal{U}(\d,\r)$ as
\begin{equation}
\mathcal{L}(\d,\r):= 1-\lmin(\delta,\rho) \quad \hbox{and}\quad \mathcal{U}(\d,\r):= \lmax(\d,\r)-1 \label{eq:LUdeltarho}
\end{equation}
with $\mathrm{H}(p):=-p\log\left(p\right)-(1-p)\log\left(1-p\right)$ denoting the usual Shannon Entropy with base $e$ logarithms, $\lmin(\d,\r)$ and $\lmax(\d,\r)$ as the solution to \eqref{eq:lmin} and \eqref{eq:lmax}, respectively:
\begin{equation}
\label{eq:lmin}
\Psi_{\min}\left(\lambda,\delta,\rho\right):=\psi_{min}(\lambda^{min}(\d,\r),\rho)+\delta^{-1}\mathrm{H}(\delta\rho)=0
\end{equation}
for  $\lambda^{min}(\d,\r)\le 1-\r$ and
\begin{equation}
\label{eq:lmax}
\Psi_{\max}\left(\lambda,\delta,\rho\right):=\psi_{max}(\lambda^{max}(\d,\r),\rho)+\delta^{-1}\mathrm{H}(\delta\rho)=0
\end{equation}
for  $\lambda^{max}(\d,\r)\ge 1+\r$ where
\begin{eqnarray}
\psi_{min}(\lambda,\rho) & := & \mathrm{H}(\rho)+\frac{1}{2}\left[(1-\rho)\log\lambda +1-\rho+\rho\log\rho-\lambda \right],\label{eq:psimin} \\
\psi_{max}(\lambda,\rho) & := & \frac{1}{2}\left[(1+\rho)\log\lambda +1+\rho-\rho\log\rho-\lambda \right].\label{eq:psimax}
\end{eqnarray}
\end{definition}

In Definition \ref{def:LUbounds}, the quantities $\psi_{min}(\lambda,\rho)$ and $\psi_{max}(\lambda,\rho)$ in \eqref{eq:psimin} and \eqref{eq:psimax}, are the large deviation exponents of the lower tail probability density function of the smallest eigenvalue and the upper tail probability density function of the largest eigenvalue of Wishart matrices respectively. The $\Psi_{\min}\left(\lambda,\delta,\rho\right)$ and $\Psi_{\max}\left(\lambda,\delta,\rho\right)$ include a Shannon entropy term from a union bound of the $\binom{N}{k}$ submatrices with $k$ columns. The level curve of $\Psi_{\min}\left(\lambda,\delta,\rho\right)$ and $\Psi_{\max}\left(\lambda,\delta,\rho\right)$ defines the transition which for $\d$ and $\r$ fixed it becomes exponentially unlikely that the smallest eigenvalue is less that $\lambda^{min}(\d,\r)$ and the largest eigenvalue is less than $\lambda^{max}(\d,\r)$.

Theorems \ref{thm:ULbctr} - \ref{thm:ULbctg} are discussed in Sections \ref{sec:fdr0} - \ref{sec:rfd} respectively.  Each section includes plots illustrating the formulae and relative difference in the relevant regimes.  The discussion of Corollary \ref{cor_ULrgd} is included in Section \ref{sec:rfd}.  This Section concludes with proofs of the compressed sensing sampling theorems discussed in Section \ref{sec:mresults}.

\subsection{Theorems \ref{thm:ULbctr}: $\delta$ fixed and $\rho\ll 1$}\label{sec:fdr0}

\begin{figure}[h]
\includegraphics*[width=0.495\textwidth]{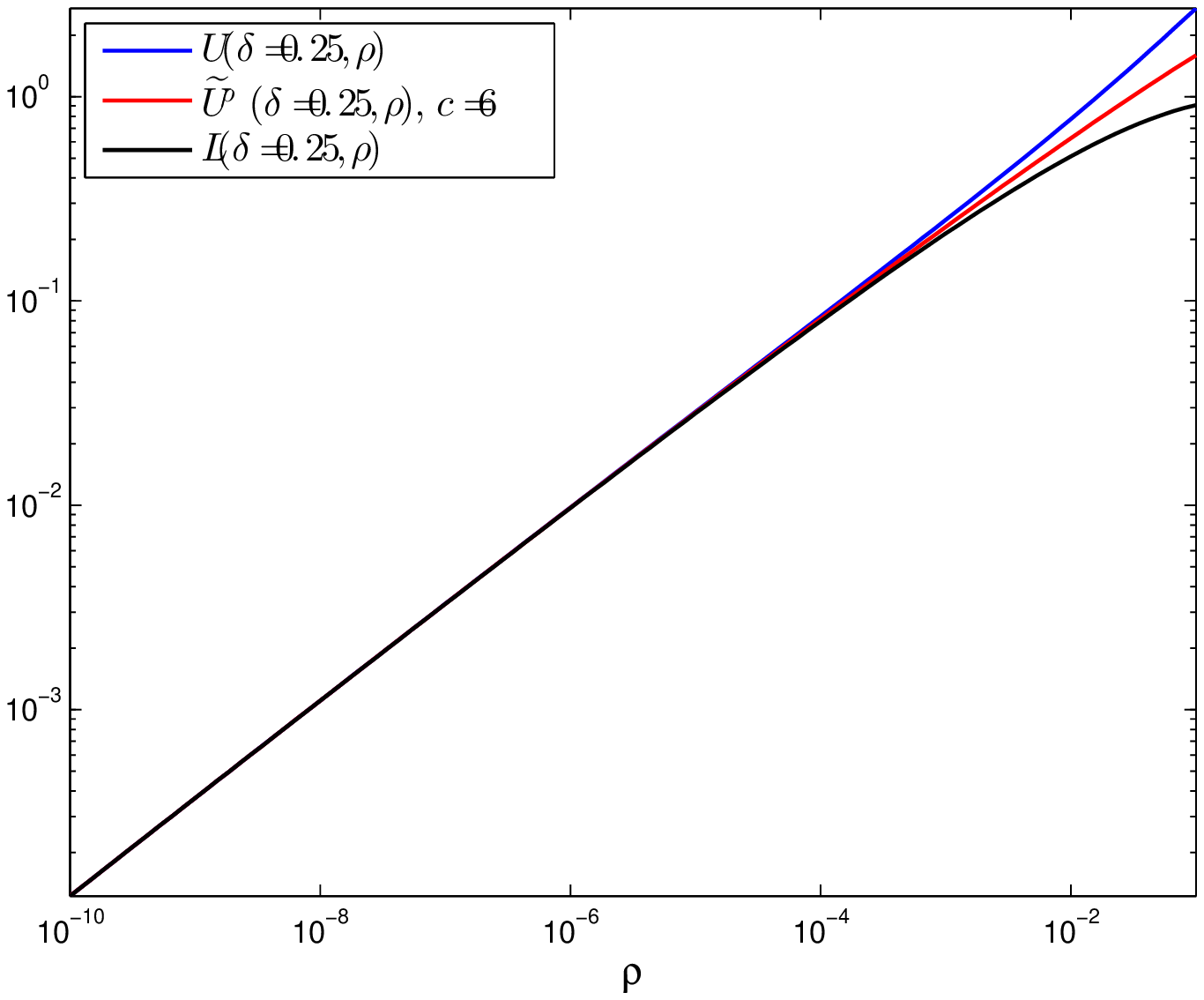}
\includegraphics*[width=0.495\textwidth]{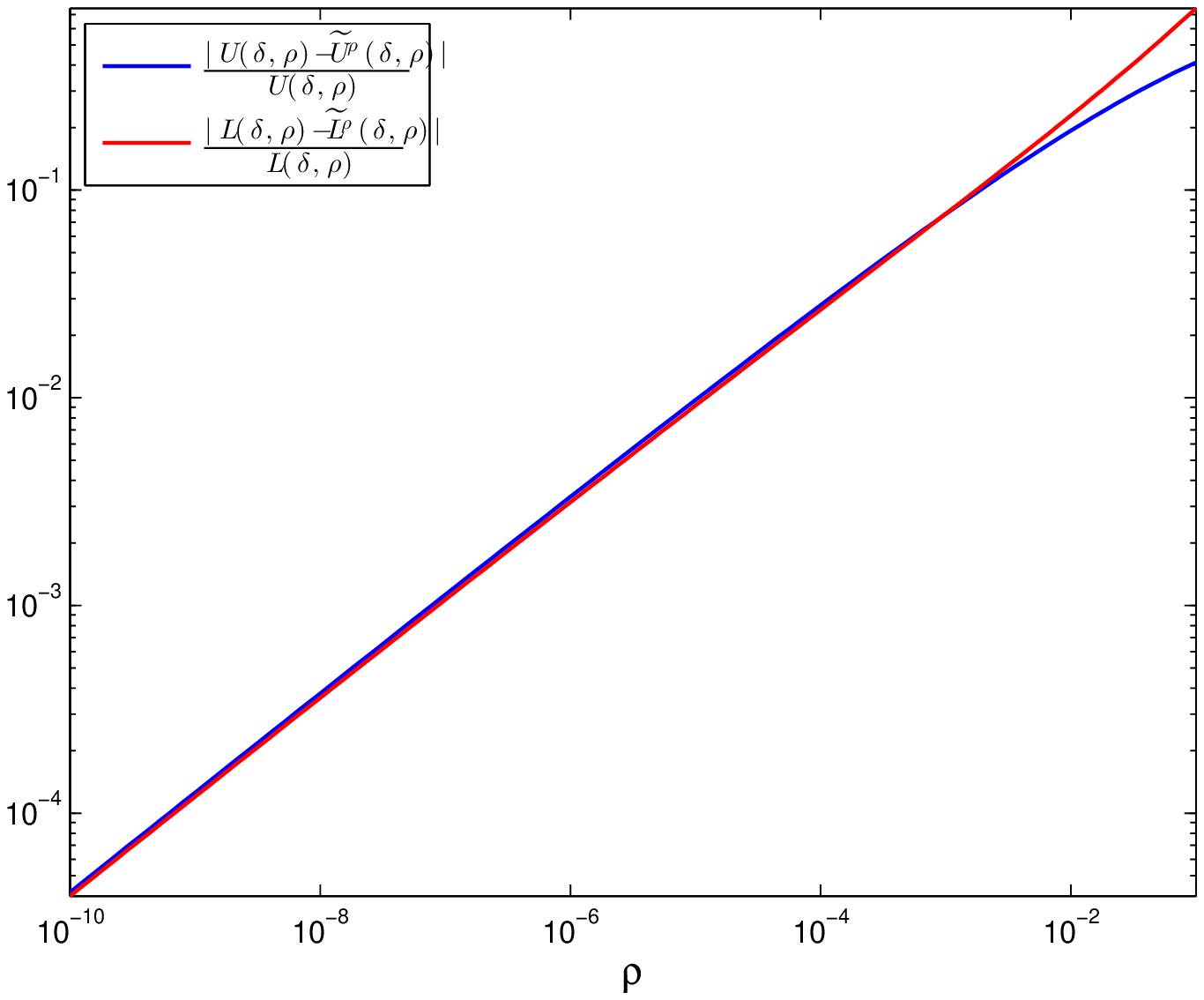}
\caption{RIC bounds for $\delta=0.25, ~c=6$ and $\rho\in(10^{-10},10^{-1})$.
{\em Left panel:} $\mathcal{U}(\delta,\rho)$, $\mathcal{L}(\delta,\rho)$, $\widetilde{\mathcal{U}}^{\r}(\delta,\rho)$ and $\widetilde{\mathcal{L}}^{\r}(\delta,\rho)$. {\em Right panel:} relative differences, $\frac{|\mathcal{U}(\delta,\rho) - \widetilde{\mathcal{U}}^{\r}(\delta,\rho)|}{\mathcal{U}(\delta,\rho)}$  and $\frac{|\mathcal{L}(\delta,\rho) - \widetilde{\mathcal{L}}^{\r}(\delta,\rho)|}{\mathcal{L}(\delta,\rho)}$.}
\label{fig1}
\end{figure}

Figure \ref{fig1}, left panel, displays the bounds $\mathcal{U}(\delta,\rho)$ and $\mathcal{L}(\delta,\rho)$ from Theorem \ref{thm:LUbounds} for $\d=0.25, ~c=6$ and $\rho\in (10^{-10},10^{-1})$.  This is the regime of Theorem \ref{thm:ULbctr} and the formulae \eqref{eq:Ubctr} and \eqref{eq:Lbctr} are also displayed.  Formulae \eqref{eq:Ubctr} and \eqref{eq:Lbctr} are observed to accurately approximate $\mathcal{U}(\delta,\rho)$ and $\mathcal{L}(\delta,\rho)$ respectively in both an absolute and relative scale, in the left and right panel of Figure \ref{fig1} respectively.

\subsection{Theorems \ref{thm:ULbctd}: $\rho$ fixed and $\delta\ll 1$}\label{sec:frd0}

\begin{figure}[h]
\centering
\includegraphics[width=0.495\textwidth]{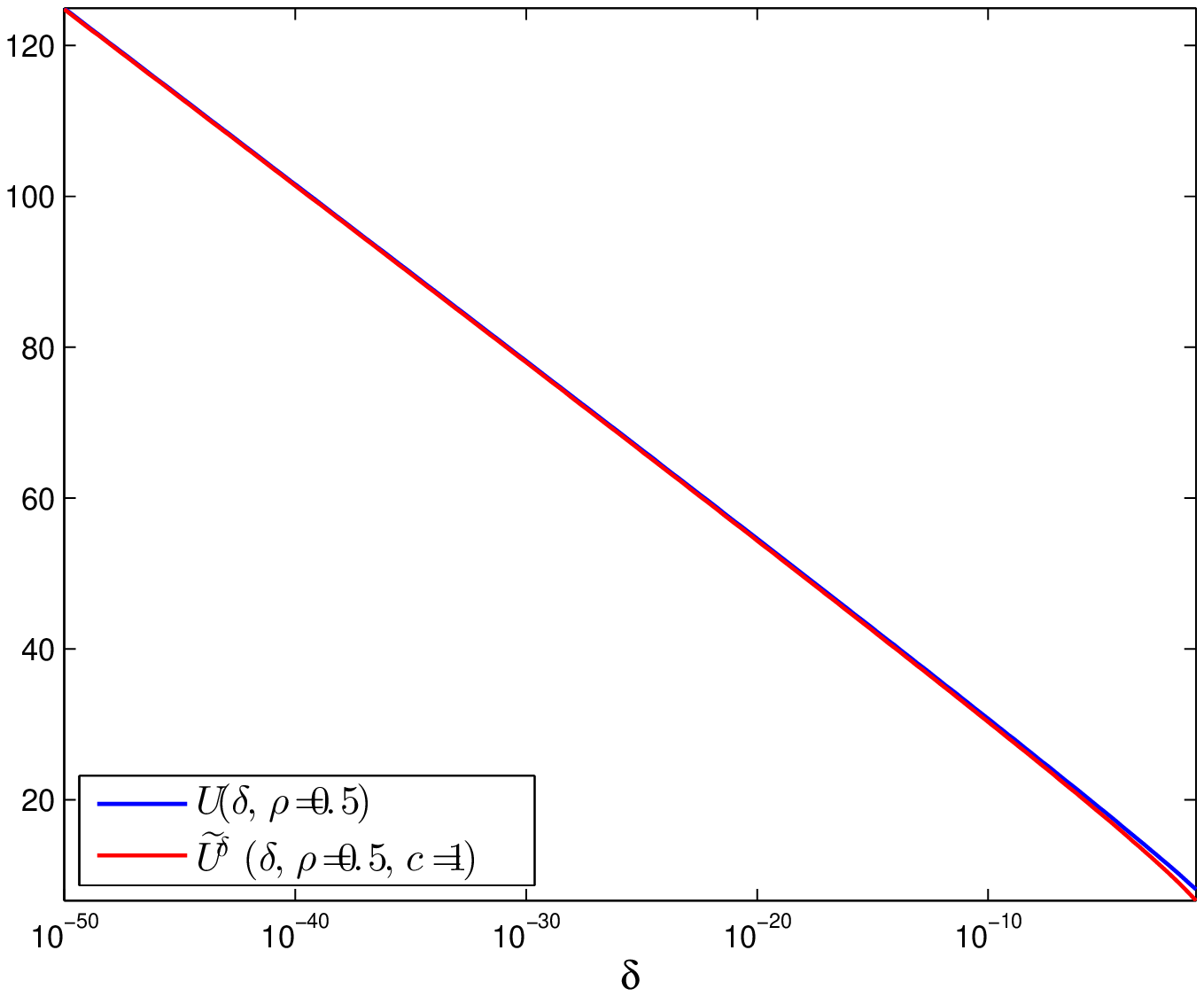}
\includegraphics[width=0.495\textwidth]{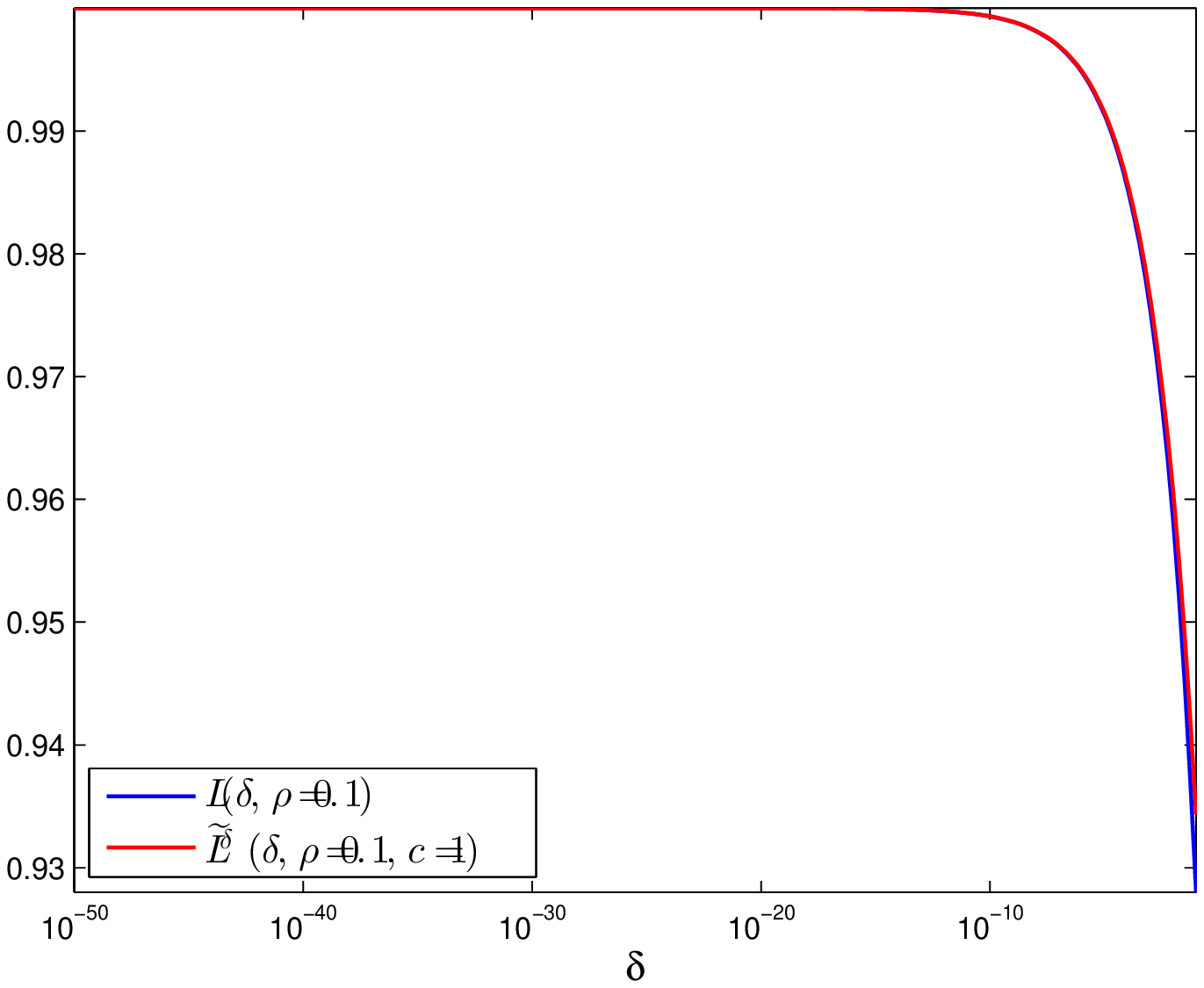}
\caption{RIC bounds for $\d\in(10^{-50},10^{-1})$ and $c=1$.  {\em Left panel:} $\mathcal{U}(\delta,\rho)$ and $\widetilde{\mathcal{U}}^{\d}(\delta,\rho)$ for $\rho = 0.5$.  {\em Right panel:} $\mathcal{L}(\delta,\rho)$ and $\widetilde{\mathcal{L}}^{\d}(\delta,\rho)$ for $\rho = 0.1$.}
\label{fig2}
\end{figure}

\begin{figure}[h]
\centering
\includegraphics[width=0.495\textwidth]{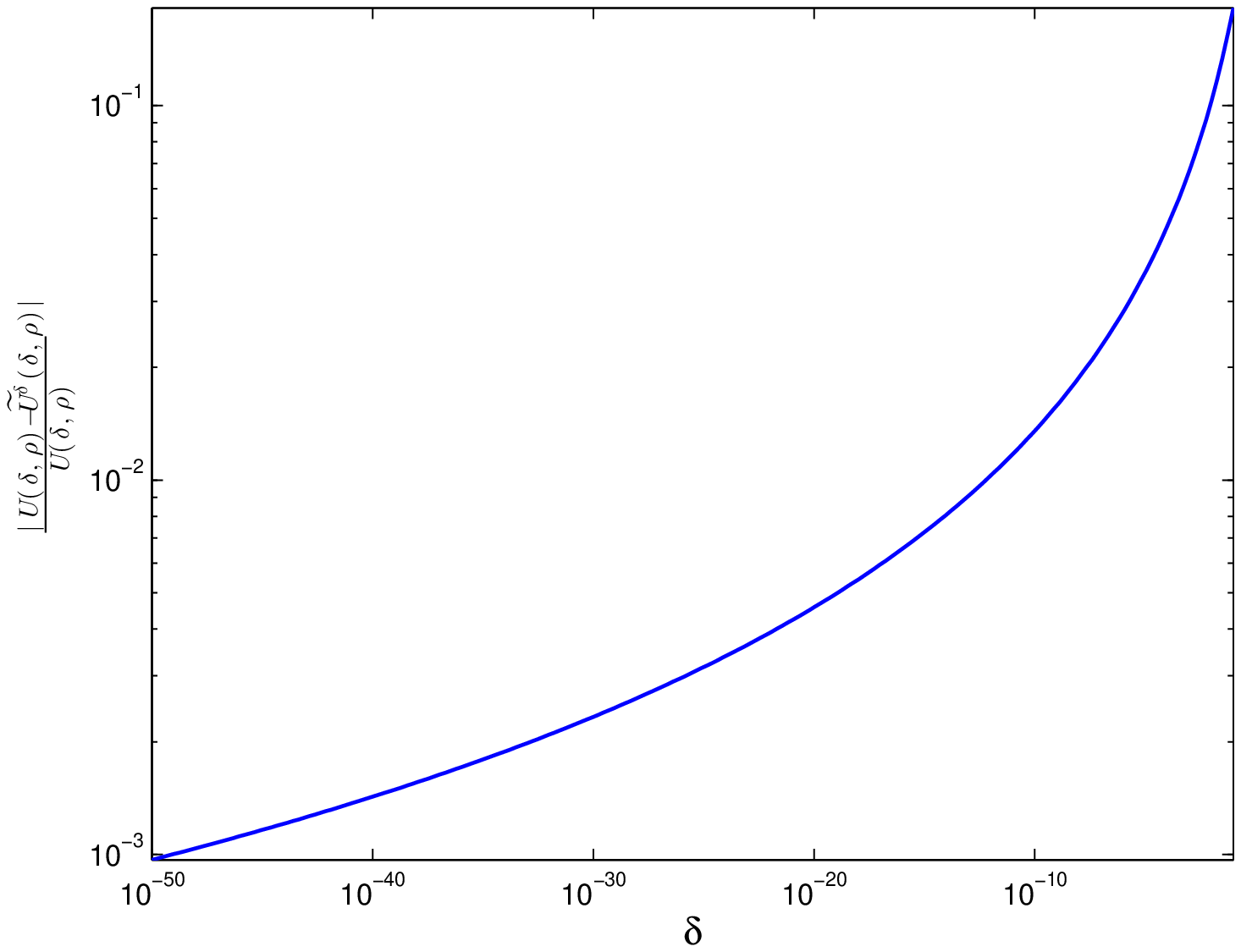}
\includegraphics[width=0.495\textwidth]{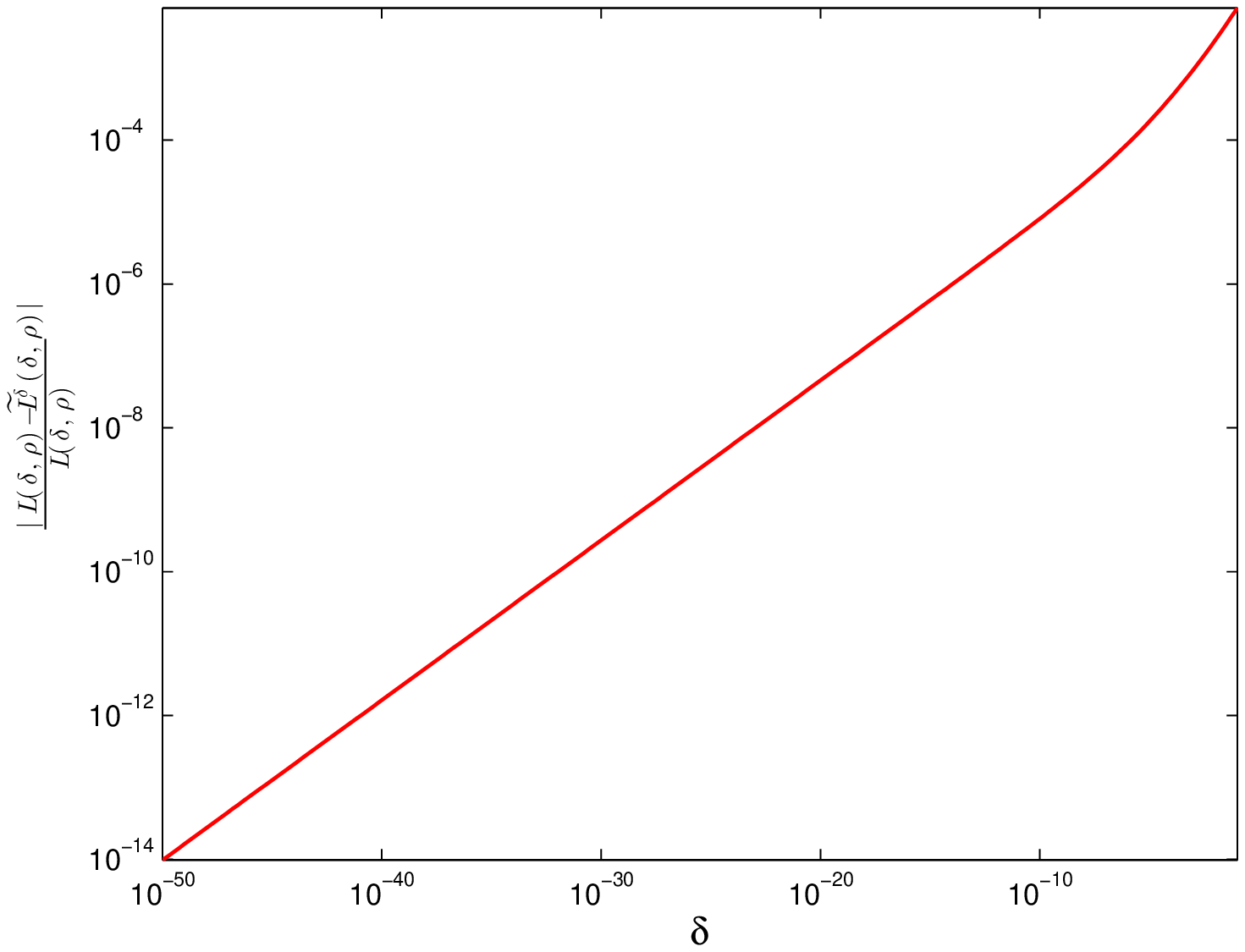}
\caption{ Relative difference in RIC bounds for $\d\in(10^{-50},10^{-1})$ and $c=1$.  {\em Left panel:} $\frac{|\mathcal{U}(\delta,\rho) - \widetilde{\mathcal{U}}^{\d}(\delta,\rho)|}{\mathcal{U}(\delta,\rho)}$ for $\rho = 0.5$. {\em Right panel:} $\frac{|\mathcal{L}(\delta,\rho) - \widetilde{\mathcal{L}}^{\d}(\delta,\rho)|}{\mathcal{L}(\delta,\rho)}$ for $\rho = 0.1$.}
\label{fig3}
\end{figure}

Figure \ref{fig2} displays the bounds $\mathcal{U}(\delta,\rho)$ and $\mathcal{L}(\delta,\rho)$ from Theorem \ref{thm:LUbounds} along with the formulae \eqref{eq:Ubctd} and \eqref{eq:Lbctd} of Theorem \ref{thm:ULbctd} in the left and right panels respectively; for diversity the upper RIC bound is shown for $\r=0.5$ and the lower RIC bound for $\r=0.1$, in both instances $\delta\in (10^{-50},10^{-1})$ and $c=1$.  This is the regime of $\rho$ fixed and $\delta\ll 1$ where the upper RIC diverges to infinity and the lower RIC converges to its trivial unit bound as $\delta$ approaches zero.  The bounds of Theorem \ref{thm:ULbctd} are observed to accurately approximate $\mathcal{U}(\delta,\rho)$ and $\mathcal{L}(\delta,\rho)$ in both an absolute and relative scale, in Figure \ref{fig2} and \ref{fig3} respectively.

\subsection{Theorems \ref{thm:ULbctg}: $\rho= (\gamma\log(1/\d))^{-1}$ and $\delta\ll 1$}\label{sec:rfd}

\begin{figure}[h]
\centering
\includegraphics[width=0.495\textwidth]{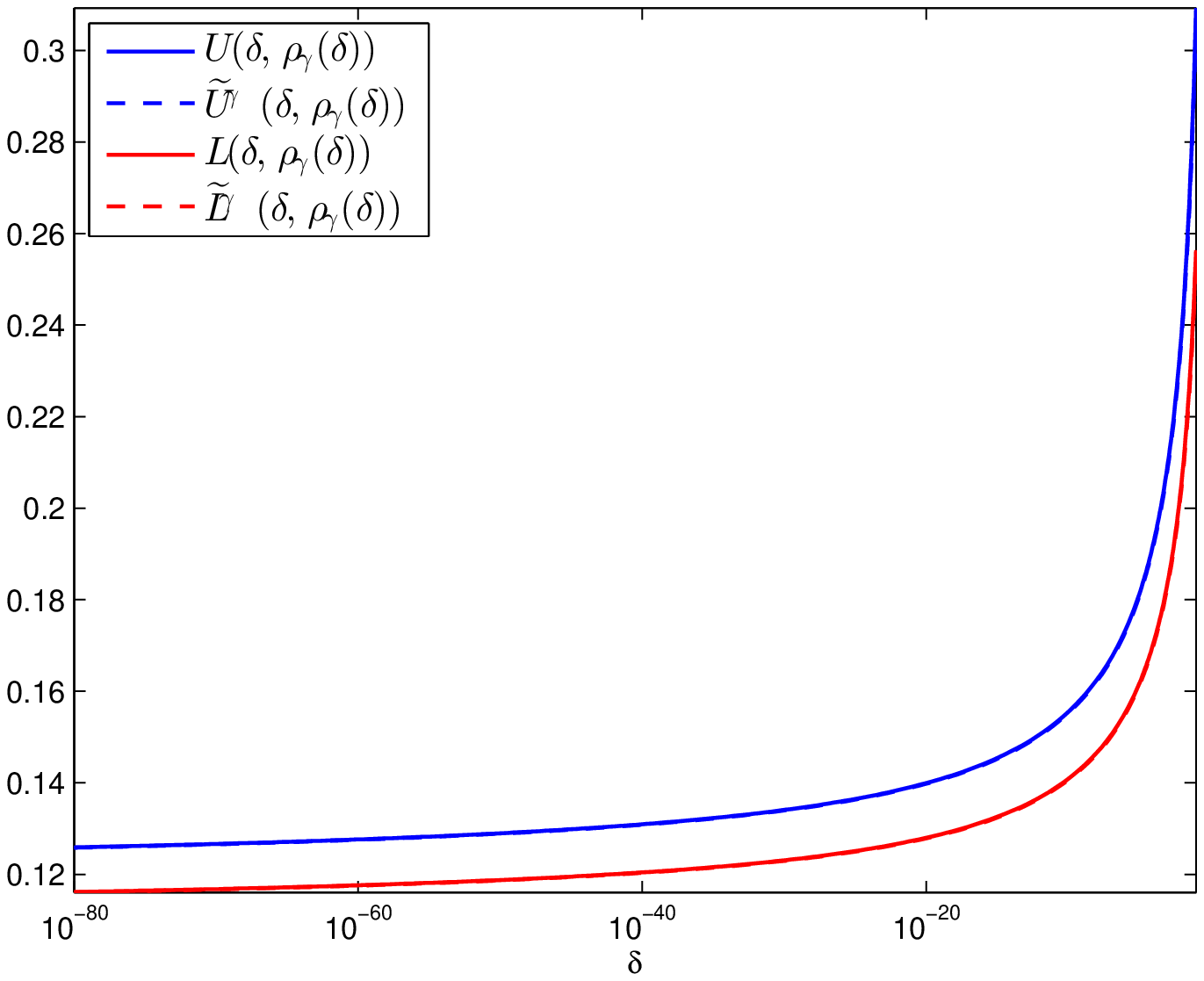}
\includegraphics[width=0.495\textwidth]{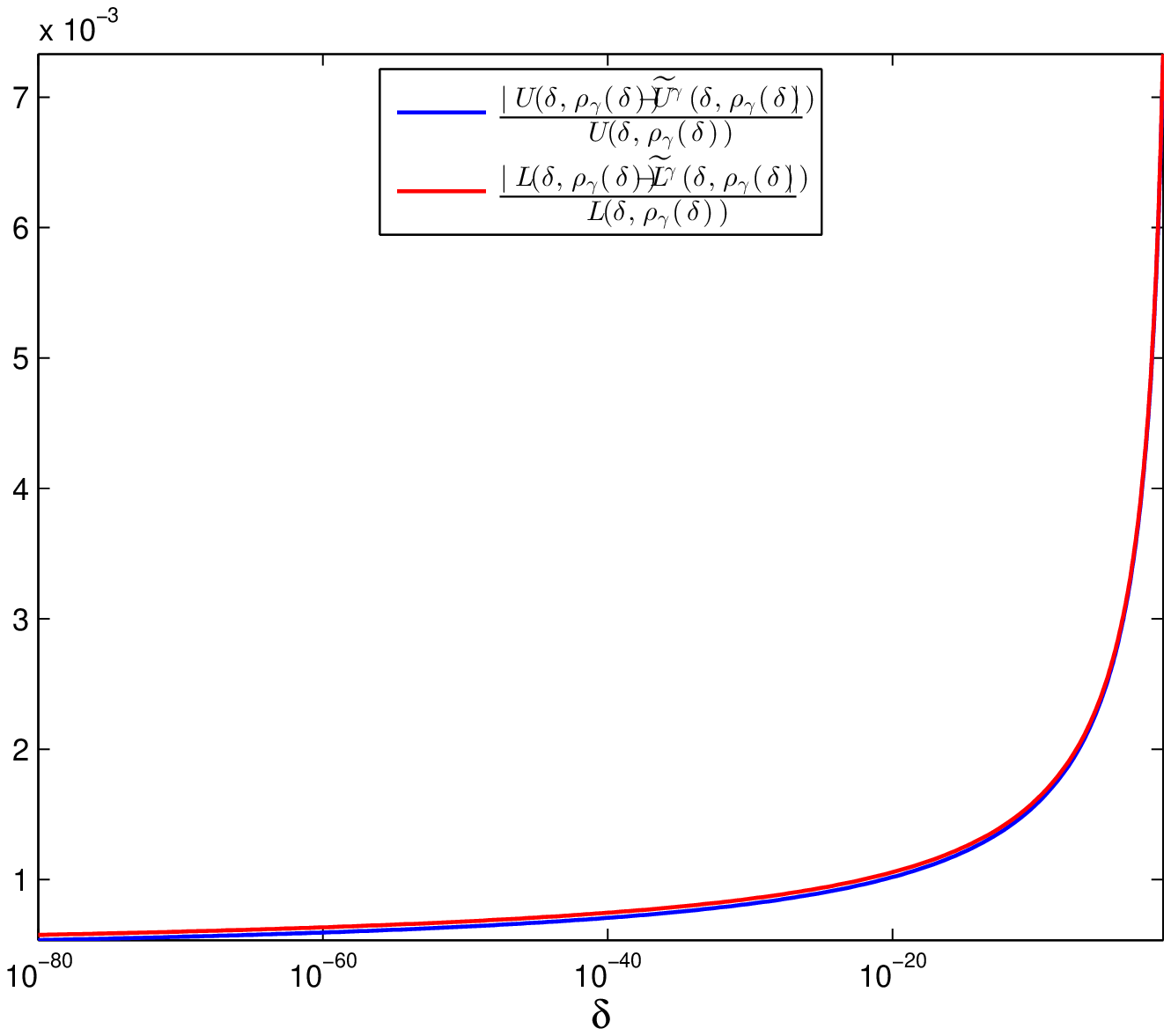}
\caption{A comparison of $\widetilde{\mathcal{U}}^{\gamma} (\delta,\rho_{\g}\left(\d)\right)$ and $\widetilde{\mathcal{L}}^{\gamma} (\delta,\rho_{\g}\left(\d)\right)$ to $\mathcal{U}(\delta,\rho)$ and $\mathcal{L}(\delta,\rho)$ respectively for $c_u = c_l = 1/3, ~\gamma = 300$ and $\delta \in (10^{-80},10^{-1})$. {\em Left panel:} $\widetilde{\mathcal{U}}^{\gamma} (\delta,\rho_{\g}\left(\d)\right)$, $\mathcal{U}(\delta,\rho)$, $\widetilde{\mathcal{L}}^{\gamma} (\delta,\rho_{\g}\left(\d)\right)$, and $\widetilde{\mathcal{L}}^{\gamma} (\delta,\rho_{\g}\left(\d)\right)$. {\em Right panel:} their relative differences $\frac{|\mathcal{U}(\delta,\rho) - \widetilde{\mathcal{U}}^{\gamma} (\delta,\rho_{\g}\left(\d)\right)|}{\mathcal{U}(\delta,\rho)}$ and $\frac{|\mathcal{L}(\delta,\rho) - \widetilde{\mathcal{L}}^{\gamma} \left(\delta,\rho_{\g}(\d)\right)|}{\mathcal{L}(\delta,\rho)}$.}
\label{fig4}
\end{figure}

The left panel of Figure \ref{fig4} displays the bounds $\mathcal{U}(\delta,\rho)$ and $\mathcal{L}(\delta,\rho)$ from Theorem \ref{thm:LUbounds} along with the formulae \eqref{eq:Ubctg} and \eqref{eq:Lbctg} of Theorem \ref{thm:ULbctg} for $c_u = c_l = 1/3, ~\gamma=300$ and  $\d\in(10^{-80},10^{-1})$.  The formulae of Theorem \ref{thm:ULbctg} are observed to accurately approximate the bounds in Theorem \ref{thm:LUbounds} over the entire range of $\delta$; the relative differences between these bounds are displayed in the right panel of Figure \ref{fig4}.

\begin{SCfigure}
\includegraphics[width=0.5\textwidth]{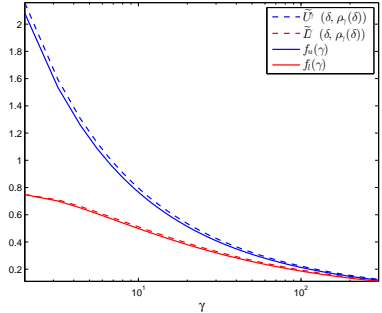}
\caption{ Plots of $\widetilde{\mathcal{U}}^{\gamma} (\delta,\rho_{\g}\left(\d)\right)$ and $\widetilde{\mathcal{L}}^{\gamma} (\delta,\rho_{\g}\left(\d)\right)$ as well as $f_u(\g)$ and $f_l(\g)$ given by \eqref{cor_Urgd} and \eqref{cor_Lrgd} respectively, for $c_u = c_l = 1/3, ~\d=10^{-80}$ and $\gamma\in (1,300)$.}
\label{fig5}
\end{SCfigure}

The left panel of Figure \ref{fig4} shows the RIC bounds converging to nonzero constants as $\d$ approaches zero, displayed for $c_u = c_l = 1/3$ and $\gamma=300$. Corollary \ref{cor_ULrgd} provides formula for $\d\ll 1$, which is observed in Figure \ref{fig5} to accurately approximate the formulae in Theorem \ref{thm:ULbctg} for $c_u = c_l = 1/3$ and $\d=10^{-80}$, uniformly over $\gamma\in(1,300)$.

\subsection{Proof of compressed sensing corollaries}\label{sec:pmresults}

Corollaries \ref{alg_corl} and \ref{omp_corl} follow directly from Theorems \ref{thm:ULbctg} and \ref{thm:ULbctr} and existing RIC based recovery guarantees for the associated algorithms in \cite{blanchard2011phase,zhang2011sparse} and \cite{mo2011remarks} respectively.

\subsubsection{Proof of Corollary \ref{alg_corl}}\label{sec:palg_corl}

\begin{proof}
There is an extensive literature on compressed sensing and sparse approximation algorithms which are guaranteed to recover vectors ${\hat x}$ that satisfy bounds of the form $\|x_0-{\hat x}\|_2\le Const\cdot \|e\|_2$ from $y=Ax_0$ provided the RICs of $A$ are sufficiently small.  The article \cite{blanchard2011phase} provides a framework by which RIC bounds can be inserted into the recovery conditions, and compressed sensing sampling theorems can be calculated from the resulting equations.  Theorem \ref{thm:ULbctg} establishes valid bounds on the RICs of Gaussian matrices in the regime considered in Corollary \ref{alg_corl}.  The claims stated in Corollary \ref{alg_corl} for $\ell_1$-minimization, IHT, SP and CoSaMP follow directly from substituting the RIC bounds of Theorem \ref{thm:ULbctg} into Theorem 10-13 of \cite{blanchard2011phase} and solving for the minimum $\gamma$ that satisfies the stated theorems. Similarly, for OMP with $31k$ steps, \cite{zhang2011sparse} provides a condition that can be expressed in the form of the framework provided by \cite{blanchard2011phase}, mentioned above. Then the claims stated in Corollary \ref{alg_corl} for OMP with $31k$ steps follows from substituting the RIC bounds of Theorem \ref{thm:ULbctg} into this condition and solving for the minimum $\gamma$. The calculated values of $\gamma$ have been rounded up to the nearest integer for ease of presentation.  Nearly identical values of $\gamma$ can be calculated using the equations from Corollary \ref{cor_ULrgd} rather than the more refined equations in Theorem \ref{thm:ULbctg}.

\end{proof}

\subsubsection{Proof of Corollary \ref{omp_corl}}\label{sec:pomp_corl}

\begin{proof}
It has been recently shown that Orthogonal Matching Pursuit (OMP) is guaranteed to recover any $k$-sparse vector after $k$ steps from its exact measurements provided, \cite{mo2011remarks},
\begin{equation}
\label{htz_omp_cond}
\max(L(k,n,N;A),U(k,n,N;A))< \frac{1}{\sqrt{k-1}}.
\end{equation}
The claimed sampling theorem is obtained by substituting the bound from Theorem \ref{thm:ULbctr} for $\max(L(k,n,N;A),U(k,n,N;A))$ and solving for $n$.
\end{proof}

\section{Proofs of Theorems \ref{thm:ULbctr} - \ref{thm:ULbctg}}\label{sec:prip_asymp}

The proof of Theorems \ref{thm:ULbctr} - \ref{thm:ULbctg} are based upon the previous analysis in \cite{bah2010improved,blanchard2011compressed}, differing in the asymptotic limits considered.  The analysis here builds upon the following large deviation bounds on the probability of the sparse eigenvalues exceeding specified values; these bounds are as follows:

With $L(k,n,N;A)$ and $U(k,n,N;A)$ defined as in \eqref{eq:L} and \eqref{eq:U} respectively, and $\Psi_{\max}\left(\lambda(\delta,\rho),\delta,\rho\right)$ and $\Psi_{\min}\left(\lambda(\delta,\rho),\delta,\rho\right)$ defined as in \eqref{eq:lmin} and \eqref{eq:lmax}, we have the bounds \cite{bah2010improved,blanchard2011compressed}

\begin{multline}
\label{ppmax}
Prob\left(\max_{K\subset\Omega,|K|=k} \lambda^{\max}(A^*_KA_K)>\lambda\right) \\ \le poly(n,\lambda)\cdot \exp\left(2n \cdot \Psi_{\max}\left(\lambda,\delta,\rho\right)\right),
\end{multline}
and
\begin{multline}
\label{ppmin}
Prob\left(\min_{K\subset\Omega,|K|=k} \lambda^{\min}(A^*_KA_K)>\lambda\right) \\ \le poly(n,\lambda)\cdot \exp\left(2n \cdot \Psi_{\min}\left(\lambda,\delta,\rho\right)\right),
\end{multline}
where $~\lambda^{\min}(B)$ and $~\lambda^{\max}(B)$ are the smallest and largest eigenvalue of $B$ respectively and $poly(z)$ is a (possibly different) polynomial function of its arguments, for explicit formulae see \cite{bah2010improved}.  Theorems \ref{thm:ULbctr} - \ref{thm:ULbctg} follow by proving that for the claimed bounds, the large deviation exponents $n\Psi_{\max}\left(\lambda(\delta,\rho),\delta,\rho\right)$ and $n\Psi_{\min}\left(\lambda(\delta,\rho),\delta,\rho\right)$ diverge to $-\infty$ as the problem size increases, and do so at a rate sufficiently fast to ensure an overall exponential decay. In addition to establishing the claims of Theorems \ref{thm:ULbctr}-\ref{thm:ULbctg}, we also show that the bounds presented in these theorems cannot be improved upon using the inequalities \eqref{ppmax} and \eqref{ppmin}, they are in fact sharp leading order asymptotic expansions of the bounds in Theorem \ref{thm:LUbounds}.

Throughout the proofs of Theorems \ref{thm:ULbctr}-\ref{thm:ULbctg} we will be using the following bounds for the Shannon entropy function, $\mathrm{H}(x):=-x\log x -(1-x)\log(1-x)$
  \begin{align}
  \label{eq:shannon_bounds}
   \mathrm{H}(x) & < -x\log x + x, \quad\mbox{and}\nonumber\\
   \mathrm{H}(x) & > -x\log x + x - x^2;
  \end{align}
 the upper bound follows from \eqref{eq:logmono_ubound} and the lower bound follows from \eqref{eq:logmm_lbound},
  \begin{align}
  \label{eq:logmono_ubound}
   -(1-x)\log(1-x) < & x \quad \forall x\in (0,1),\\
  \label{eq:logmm_lbound}
   -\log(1-x) > & x \quad \forall x<1 \quad \mbox{and} \quad x \neq 0.
  \end{align}

\subsection{Theorem \ref{thm:ULbctr}}\label{sec:pthm:ULbctr}

\subsubsection{The upper bound, $\widetilde{\mathcal{U}}^{\r}(\delta,\rho)$}\label{sec:pUbctrthm}

\begin{proof}

Define $$\widetilde{\lambda}^{\max}_{\r}(\delta,\rho) := 1 + \sqrt{2\rho\log\left(\frac{1}{\delta^2\rho^3}\right) + c\rho}, \quad \Rightarrow \quad \widetilde{\mathcal{U}}^{\r}(\delta,\rho) = \widetilde{\lambda}^{\max}_{\r}(\delta,\rho) - 1$$ as from \eqref{eq:Ubctr}. Bounding $\widetilde{\mathcal{U}}^{\r}(\delta,\rho)$ from above by $(1+\e)\widetilde{\mathcal{U}}^{\r}(\delta,\rho)$ is equivalent to bounding from above $\widetilde{\lambda}^{\max}_{\r}$ by $(1+\e)\widetilde{\lambda}^{\max}_{\r}-\epsilon$. We first establish that for a slightly looser bound, with $c>6$, the exponent $\Psi_{\max}\left((1+\e)\widetilde{\lambda}^{\max}_{\r}-\epsilon,\delta,\rho\right)$ is negative, and then verify that when multiplied by $n$ it diverges to $-\infty$ as $n$ increases. We also show that for a slightly tighter bound, with $c<6$, $\Psi_{\max}\left((1-\e)\widetilde{\lambda}^{\max}_{\r}+\epsilon,\delta,\rho\right)$ is positive, and hence the bound $\widetilde{\mathcal{U}}^{\r}(\delta,\rho)$ cannot be improved using the inequality \eqref{ppmax} from \cite{blanchard2011compressed}.  We show the above properties, in two parts that for $\delta$ fixed:
\begin{enumerate}
  \item $\exists ~\rho_0, ~\epsilon > 0 ~\& ~c>6$ such that for $\rho < \rho_0, \Psi_{\max}\left((1+\e)\widetilde{\lambda}^{\max}_{\r}-\epsilon,\delta,\rho\right) \leq 0;$
  \item $\nexists ~\rho_0, ~\epsilon > 0 ~\& ~c<6$ such that for $\rho < \rho_0, \Psi_{\max}\left((1-\e)\widetilde{\lambda}^{\max}_{\r}+\epsilon,\delta,\rho\right) \leq 0,$
\end{enumerate}
which are proven below separately as Part 1 and Part 2 respectively.

\begin{description}
  \item[Part 1:]
  \begin{multline}
  \label{part1r_1}
  2\Psi_{\max}\left((1+\e)\widetilde{\lambda}^{\max}_{\r}-\epsilon,\delta,\rho\right) = (1+\rho)\log\left((1+\e)\widetilde{\lambda}^{\max}_{\r}-\epsilon\right) \\
  - \rho\log(\rho) + \rho + 1 - \left((1+\e)\widetilde{\lambda}^{\max}_{\r}-\epsilon\right)+ \frac{2}{\delta}\mathrm{H}(\delta\rho),	
  \end{multline}
  by substituting $(1+\e)\widetilde{\lambda}^{\max}_{\r}-\epsilon$ for $\lambda$ in \eqref{eq:lmax}. We consolidate notation using $u := \widetilde{\lambda}^{\max}_{\r} - 1$ and
 using the first bounds of the Shannon entropy in \eqref{eq:shannon_bounds} we bound \eqref{part1r_1} above as follows
  \begin{align}
  &   2\Psi_{\max}\left((1+\e)\widetilde{\lambda}^{\max}_{\r}-\epsilon,\delta,\rho\right) \nonumber\\
  \label{part1r_2}
  & < (1+\rho)\log\left[(1+\e)(1+u) - \e\right] - \rho\log\r + \rho + 1 - (1+\e)(1+u)\nonumber\\
  & \quad + \e + \frac{2}{\d}\left[-\d\r\log\left(\d\r\right) + \d\r\right],\\
  \label{part1r_3}
  & = (1+\rho)\log\left[1 + (1+\e)u \right] + \rho\log\left(\frac{1}{\d^2\r^3}\right) + \rho - u - \e u + 2\d\r.
  \end{align}
  From \eqref{part1r_2} to \eqref{part1r_3} we expanded the products of $(1+\e)(1+u)$ and simplified.

  Now replacing $\rho\log\left(\frac{1}{\delta^2\rho^3}\right)$ by its equivalent $\left(u^2-c\r\right)/2$ and expanding $(1+\r)$ in the first term we bound \eqref{part1r_3} by
  \begin{align}
  & 2\Psi_{\max}\left((1+\e)\widetilde{\lambda}^{\max}_{\r}-\epsilon,\delta,\rho\right) \nonumber \\
  \label{part1r_4}
  & < \log\left(1 + u +\e u \right) + \rho\log\left(1 + u +\e u \right) + \frac{1}{2}\left(u^2-c\r\right) + 3\r - u - \e u, \\
  \label{part1r_5}
  & = \log(1 + u) + \log\left(1 + \frac{\e u}{1+u} \right) + \frac{1}{2}u^2 - \frac{1}{2}c\r + 3\r - u - \e u \nonumber\\
  & \quad + \rho\log\left(1 + u \right) + \rho\log\left(1 + \frac{\e u}{1+u} \right), \\
  \label{part1r_6}
  & < u - \frac{1}{2}u^2 + \frac{1}{3}u^3 + \frac{\e u}{1+u} + \frac{1}{2}u^2 - \frac{1}{2}(c-6) \r - u - \e u\nonumber\\
  & \quad + \rho u + \frac{\e\rho u}{1+u}.
  \end{align}

  From \eqref{part1r_4} to \eqref{part1r_5} the term $\log(1+u+\e u)$ is factored as in the first two logarithms in \eqref{part1r_5}.  From \eqref{part1r_5} to \eqref{part1r_6} we bounded the first  $\log(1+u)$ from above using the second bound in \eqref{eq:logp_ubound} and bounded above all other logarithmic terms using the first bound in \eqref{eq:logp_ubound} .

\begin{eqnarray}
\log(1+x) & \leq  & x, \label{eq:logp_ubound} \\
\log(1+x) & \leq  &x- \frac{1}{2}x^2 + \frac{1}{3}x^3  \quad \forall x> -1.\nonumber
% \\ \log(1+x) & \leq  &x- \frac{1}{2}x^2 + \frac{1}{3}x^3- \frac{1}{2}x^4 + \frac{1}{3}x^5 \quad \forall x> -1. \nonumber \qquad
  \end{eqnarray}

 We can bound above $1/(1+u)$ in the fourth and last terms of \eqref{part1r_6} using the bound of \eqref{eq:r_pfrac_series} below.
  \begin{align}
  \label{eq:r_pfrac_series}
  \frac{1}{1+x} & < 1 \quad \mbox{for} \quad 0<x<1.
% ,\\   \frac{1}{1+x} & < 1 - x + x^2 \nonumber
  \end{align}

Therefore, \eqref{part1r_6} becomes
 \begin{align}
& 2\Psi_{\max}\left((1+\e)\widetilde{\lambda}^{\max}_{\r}-\epsilon,\delta,\rho\right) \nonumber \\
\label{part1r_7}
  & < \frac{1}{3}u^3 - \frac{1}{2}(c-6)\r - \e u + \e u + \rho u + \e\rho u, \qquad\\
   \label{part1r_8}
  & = - \frac{1}{2}(c-6)\r + \frac{1}{3}u^3 + (1 + \e)\rho u, \\
   \label{part1r_9}
  & < - \frac{1}{4}(c-6)\r - \frac{1}{4}(c-6)\r + \frac{1}{3}u^3 + \frac{1}{14}(1+\e)u^3, \\
   \label{part1r_10}
  & = - \frac{1}{4}(c-6)\r - \frac{1}{4}(c-6)\r + \frac{17+3\e}{42}u^3.
  \end{align}

  We simplified \eqref{part1r_7} to get \eqref{part1r_8}. From \eqref{part1r_8} to \eqref{part1r_9} we split the first term into half and bounded above $\r u$ by $u^2/14$ using the fact that by the definition of $u$,
  \begin{equation*}
   u^2 = \r\left[2\log\left(\frac{1}{\d^2\r^3}\right) + 7\right] \quad \Rightarrow \quad \frac{1}{4\log\left(\frac{1}{\d^2\r^3}\right)}u^2<\r< \frac{1}{14}u^2.
  \end{equation*}
  Then we simplified from \eqref{part1r_9} to \eqref{part1r_10}.

  Now in \eqref{part1r_10}, if the sum of the last two terms is non-positive there would be a unique $\r_0$ such that as $\r \rightarrow 0$ for any $\r<\r_0$ and fixed $\d$ \eqref{part1r_10} will be negative. This is achieved if $c>6$ and
 \begin{equation}
  \label{eq:fr_p1r_1}
   - \frac{1}{4}(c-6)\r + \frac{17+3\e}{42}u^3 \leq 0 \quad \Rightarrow \quad u^3 \leq \frac{21(c-6)}{2(17+3\e)}\r.
  \end{equation}

  Since $u$ is strictly decreasing in $\rho$, there is a unique $\rho_0$ that satisfies \eqref{eq:fr_p1r_1} and makes \eqref{part1r_10} negative for $\d$ fixed, $\e>0, ~c>6$ and $\r<\r_0$ as $\r \rightarrow 0$.

  Having established a negative bound from above and the $\r_0$ for which it is valid, it remains to show that $~n \cdot 2\Psi_{\max}\left((1+\e)\widetilde{\lambda}^{\max}_{\r}-\epsilon,\delta,\rho\right) \rightarrow -\infty$ as $(k,n,N) \rightarrow \infty$.  The claimed exponential decay with $k$ follows by noting that $n\cdot\r = k$, which in conjunction with the first term in the right hand side of \eqref{part1r_10} gives a concluding bound $-(c-6) k/4.$ For $\rho<\rho_0$ therefore
  \begin{multline*}
  Prob\left(U(k,n,N;A) > (1+\e)\widetilde{\mathcal{U}}^{\r}(\delta,\rho)\right) \\ \le poly\left(n,(1+\e)\widetilde{\lambda}^{\max}_{\r}-\epsilon\right)\cdot\exp\left[-\frac{(c-6)k}{4}\right].
  \end{multline*}
  The above bound goes to zero as $k\rightarrow \infty$ provided $(\log n)/k\rightarrow 0$ so that the exponential decay in $k$ dominates the polynomial decrease in $n$.\\

  \item[Part 2:]
  \begin{multline}
  \label{part2r_1}
  2\Psi_{\max}\left((1-\e)\widetilde{\lambda}^{\max}_{\r}+\epsilon,\delta,\rho\right) = (1+\rho)\log\left((1-\e)\widetilde{\lambda}^{\max}_{\r}+\epsilon\right) \\
  - \rho\log(\rho) + \rho + 1 - \left((1-\e)\widetilde{\lambda}^{\max}_{\r}+\epsilon\right) + \frac{2}{\delta}\mathrm{H}(\delta\rho),	
  \end{multline}
  by substituting $(1-\e)\widetilde{\lambda}^{\max}_{\r}+\epsilon$ for $\lambda$ in \eqref{eq:lmax}. We consolidate notation using $u:= \widetilde{\lambda}^{\max}_{\r} - 1$ and bound the Shannon entropy function from below using the second bound in \eqref{eq:shannon_bounds} to give
  \begin{align}
  &   2\Psi_{\max}\left((1-\e)\widetilde{\lambda}^{\max}_{\r}+\epsilon,\delta,\rho\right) \nonumber\\
  \label{part2r_2}
  & > (1+\rho)\log\left[(1-\e)(1+u) + \e\right] - \rho\log\r + \rho + 1 - (1-\e)(1+u)\nonumber\\
  & \quad - \e + \frac{2}{\d}\left[-\d\r\log\left(\d\r\right) + \d\r - \d^2\r^2\right],\\
  \label{part2r_3}
  & = (1+\rho)\log\left[1 + (1-\e)u\right] + \rho\log\left(\frac{1}{\d^2\r^3}\right) + 3\rho - (1-\e) u- 2\d\r^2.
  \end{align}
  From \eqref{part2r_2} to \eqref{part2r_3} we expanded the products of $(1-\e)(1+u)$ and simplified.

  Now replacing $\rho\log\left(\frac{1}{\delta^2\rho^3}\right)$ by $\left(u^2-c\r\right)/2$ and expanding $(1+\r)$ in the first term we have \eqref{part2r_3} become
  \begin{align}
  & 2\Psi_{\max}\left((1-\e)\widetilde{\lambda}^{\max}_{\r}+\epsilon,\delta,\rho\right) \nonumber \\
  \label{part2r_4}
  & > \log\left[1 + (1-\e)u\right] + \rho\log\left[1 + (1-\e)u\right] + \frac{1}{2}\left(u^2-c\r\right) + 3\r\nonumber\qquad\\
  & \quad - (1-\e) u- 2\d\r^2, \\
  \label{part2r_5}
   & > (1-\e)u- \frac{(1-\e)^2}{2} u^2 + \frac{1}{2}u^2 - \frac{1}{2}c\r + 3\r - (1-\e) u+ \r(1-\e)u\nonumber\\
  & \quad - \frac{(1-\e)^2}{2}\r u^2 - 2\d\r^2,\\
  \label{part2r_6}
   & = \frac{\e(2-\e)}{2} u^2 + \frac{1}{2}(6-c)\r + \r u- \e\r u- \frac{(1-\e)}{2}\r u - 2\d\r^2,\qquad \\
   \label{part2r_7}
   & > \frac{1}{2}(6-c)\r + \frac{1-\e}{2}\r u- 2\d \r^2.
  \end{align}
  From \eqref{part2r_4} to \eqref{part2r_5} we bounded below the logarithmic terms by the first two terms of their series expansion using \eqref{eq:logp_lbound}
  \begin{equation}
  \label{eq:logp_lbound}
  \log(1+x) \geq  x - \frac{1}{2}x^2 \quad \forall x > -1. \qquad
  \end{equation}
 From \eqref{part2r_5} to \eqref{part2r_6} we bounded above $\r u^2$ and $(1-\e)^2$ by $\r u$ and $1-\e$ respectively and simplified. Then we dropped the first term to bound below \eqref{part2r_6} by \eqref{part2r_7} and we simplified the terms with $\r u$.

  For $c<6$, the only negative term in \eqref{part2r_7}, the last term, goes faster to zero than the rest. Therefore, there does not exist a $\r_0,~\e>0$ and $c<6$ such that for $\r <\r_0$ and fixed $\d$ \eqref{part2r_7} is negative. Thus the bound
  \begin{multline*}
  Prob\left(U(k,n,N;A) > (1-\e)\widetilde{\mathcal{U}}^{\r}(\delta,\rho)\right) \\ \le poly\left(n,(1-\e)\widetilde{\lambda}^{\max}_{\r}+\epsilon\right)\cdot\exp\left[2n\Psi_{\max}\left((1-\e)\widetilde{\lambda}^{\max}_{\r}+\epsilon,\delta,\rho\right)\right],
  \end{multline*}
  does not decay to zero as $n\rightarrow\infty$.

\end{description}

Now \textbf{Part 1} and \textbf{Part 2} put together shows that $\widetilde{\mathcal{U}}^{\r}(\delta,\rho)$ is a tight upper bound of $U(k,n,N;A)$ with overwhelming probability as the problem size grows in the regime prescribed for $\widetilde{\mathcal{U}}^{\r}(\delta,\rho)$ in Theorem \ref{thm:ULbctr}.

\end{proof}

\subsubsection{The lower bound, $\widetilde{\mathcal{L}}^{\r}(\delta,\rho)$}\label{sec:pLbctrthm}

\begin{proof}
Define $$\widetilde{\lambda}^{\min}_{\r}(\delta,\rho) := 1 - \sqrt{2\rho\log\left(\frac{1}{\delta^2\rho^3}\right) + c\rho},\quad \Rightarrow \quad \widetilde{\mathcal{L}}^{\r}(\delta,\rho) = 1 - \widetilde{\lambda}^{\min}_{\r}(\delta,\rho)$$ as from \eqref{eq:Lbctr}. Since bounding $\widetilde{\mathcal{L}}^{\r}(\delta,\rho)$ above by $(1+\e)\widetilde{\mathcal{L}}^{\r}(\delta,\rho)$ is equivalent to bounding $\widetilde{\lambda}^{\min}_{\r}$ above by $(1+\e)\widetilde{\lambda}^{\min}_{\r}-\epsilon$. We first establish that for a slightly looser bound, with $c>6$, the exponent $\Psi_{\min}\left((1+\e)\widetilde{\lambda}^{\min}_{\r}-\epsilon,\delta,\rho\right)$, and then verify that when multiplied by $n$ it diverges to $-\infty$ as $n$ increases.
We also show that for a slightly tighter bound, with $c<6$, $\Psi_{\min}\left((1-\e)\widetilde{\lambda}^{\min}_{\r}+\epsilon,\delta,\rho\right)$ is positive, and hence the bound $\widetilde{\mathcal{L}}^{\r}(\delta,\rho)$ cannot be improved using the inequality \eqref{ppmin} from \cite{blanchard2011compressed}.  We show, in two parts that for $\delta$ fixed:
\begin{enumerate}
  \item $\exists ~\rho_0, ~\epsilon > 0 ~\& ~c>6$ such that for $\rho < \rho_0, \Psi_{\min}\left((1+\e)\widetilde{\lambda}^{\min}_{\r}-\epsilon,\delta,\rho\right) \leq 0;$
  \item $\nexists ~\rho_0, ~\epsilon > 0 ~\& ~c<6$ such that for $\rho < \rho_0, \Psi_{\min}\left((1-\e)\widetilde{\lambda}^{\min}_{\r}+\epsilon,\delta,\rho\right) \leq 0,$
\end{enumerate}
which are proven separately in the two parts as follows.

\begin{description}
  \item[Part 1:]
  \begin{multline}
  \label{part1r_1b}
  2\Psi_{\min}\left((1+\e)\widetilde{\lambda}^{\min}_{\r}-\epsilon,\delta,\rho\right) = 2\mathrm{H}(\rho) + (1-\rho)\log\left((1+\e)\widetilde{\lambda}^{\min}_{\r}-\epsilon\right)\\ + \rho\log(\rho) - \rho + 1 - \left((1+\e)\widetilde{\lambda}^{\min}_{\r}-\epsilon\right)+ \frac{2}{\delta}\mathrm{H}(\delta\rho),	
  \end{multline}
  by substituting $(1+\e)\widetilde{\lambda}^{\min}_{\r}-\epsilon$ for $\lambda$ in \eqref{eq:lmin}. We consolidate notation using $l := 1 - \widetilde{\lambda}^{\min}_{\r}$ and bound the Shannon entropy functions from above using the first bound in \eqref{eq:shannon_bounds} which gives
  \begin{align}
  & 2\Psi_{\min}\left((1+\e)\widetilde{\lambda}^{\min}_{\r}-\epsilon,\delta,\rho\right) \nonumber\\
  \label{part1r_2b}
  & < -2\r\log\left(\r\right) + 2\r + (1-\rho)\log\left[(1+\e)(1-l) - \e\right] + \rho\log\r \nonumber\\
  & \quad - \rho + 1 - (1+\e)(1-l) + \e - 2\r\log\left(\d\r\right) + \frac{2}{\d}(\d\r) ,\\
  \label{part1r_3b}
  & = (1-\rho)\log\left(1-l-\e l \right) + \rho\log\left(\frac{1}{\d^2\r^3}\right) + 3\rho + l + \e l.
  \end{align}
  We simplified from \eqref{part1r_2b} to \eqref{part1r_3b}.

  Now replacing $\rho\log\left(\frac{1}{\delta^2\rho^3}\right)$ by $\left(l^2-c\r\right)/2$ and factoring $(1-l)$ in the argument of the first log term we have \eqref{part1r_3b} become
  \begin{align}
  & 2\Psi_{\min}\left((1+\e)\widetilde{\lambda}^{\min}_{\r}-\epsilon,\delta,\rho\right) \nonumber \\
  \label{part1r_4b}
  & < (1-\r)\log(1 - l) + (1-\rho)\log\left(1 - \frac{\e l}{1-l} \right) + \frac{1}{2}\left(l^2-c\r\right) + 3\r\nonumber\\
  & \quad + l + \e l, \\
  \label{part1r_5b}
  & < l + \log(1 - l) + \frac{1}{2}l^2 - \frac{1}{2}c\r + 3\r -\r\log(1 - l) + \e l\nonumber\\
  & \quad - (1-\rho)\frac{\e l}{1-l}, \\
  \label{part1r_6b}
  & < l - l - \frac{1}{2}l^2 + \frac{1}{2}l^2 - \frac{1}{2}(c-6)\r -\r\log(1 - l) + \e l - \e l (1-\rho), \\
  \label{part1r_7b}
  & = - \frac{1}{2}(c-6)\r -\r\log(1 - l) + \e l - \e l + \e\r l, \\
  \label{part1r_8b}
  & = - \frac{1}{4}(c-6)\r - \frac{1}{4}(c-6)\r -\r\log(1 - l) + \e\rho l.
  \end{align}
  From \eqref{part1r_4b} to \eqref{part1r_5b} we expanded $(1-\r)$ and we bounded above the second logarithmic term using the first bound of \eqref{eq:logm_ubound}.
  \begin{eqnarray}
   \log(1-x) & \leq & -x, \label{eq:logm_ubound}\\
   \log(1-x) & \leq & -x - \frac{1}{2}x^2, \nonumber\\
   \log(1-x) & \leq & -x - \frac{1}{2}x^2 - \frac{1}{3}x^3 \quad \forall x\in(0,1). \nonumber
  \end{eqnarray}
  From \eqref{part1r_5b} to \eqref{part1r_6b} we bounded above the first logarithmic term using the second bound of \eqref{eq:logm_ubound} and also bounded $1/(1-l)$ using \eqref{eq:r_mfrac_series}.
  \begin{equation}
   \frac{1}{1-x} \geq 1 \quad \forall x \in (0,1). \label{eq:r_mfrac_series}
  \end{equation}

  From \eqref{part1r_6b} to \eqref{part1r_7b} we expanded the last brackets and simplified and from \eqref{part1r_7b} to \eqref{part1r_8b} we simplified and split the first term into two equal terms.

  Equation \eqref{part1r_8b} is clearly negative if $c>6$ and the sum of the last three terms is non-positive, which is satisfied if $\epsilon l - \log(1- l)\le (c-6)/4$, which is also true if, using the first bound in \eqref{eq:logp_ubound}, $(1+\epsilon) l\le (c-6)/4$.  Since $l$ is strictly increasing in $\rho$, taking on values between zero and 1, there is a unique $\rho_0$ such that for fixed $\d$, $\e>0$ and $c>6$, any $\r<\r_0$ satisfies $(1+\epsilon) l\le (c-6)/4$ and \eqref{part1r_8b} is negative.

  Having established a negative bound from above and the $\r_0$ for which it is valid, it remains to show that $~n \cdot 2\Psi_{\min}\left((1+\e)\widetilde{\lambda}^{\min}_{\r}-\epsilon,\delta,\rho\right) \rightarrow -\infty$ as $(k,n,N) \rightarrow \infty$, which verifies an exponential decay to zero of the bound \eqref{ppmin} with $k$. This follows by noting that $n\cdot\r = k$, which in conjunction with the first term in the right hand side of \eqref{part1r_8b} gives a concluding bound $-\frac{1}{4}(c-6) k.$ For $\rho<\rho_0$ therefore
  \begin{multline*}
  Prob\left(L(k,n,N;A) > (1+\e)\widetilde{\mathcal{L}}^{\r}(\delta,\rho)\right) \\ \le poly\left(n,(1+\e)\widetilde{\lambda}^{\min}_{\r}-\epsilon\right)\cdot\exp\left[-\frac{(c-6)k}{4}\right].
  \end{multline*}
  The right hand side of which goes to zero as $k\rightarrow \infty$ with $(\log n)/k\rightarrow 0$ as $k\rightarrow \infty$ so that the exponential decay in $k$ dominates the polynomial decrease in $n$.\\

  \item[Part 2:]
  \begin{multline}
  \label{part2r_1b}
  2\Psi_{\min}\left((1-\e)\widetilde{\lambda}^{\min}_{\r}+\epsilon,\delta,\rho\right) = 2\mathrm{H}(\rho) + (1-\rho)\log\left((1-\e)\widetilde{\lambda}^{\min}_{\r}+\epsilon\right)\\ + \rho\log(\rho) - \rho + 1 - \left((1-\e)\widetilde{\lambda}^{\min}_{\r}+\epsilon\right)+ \frac{2}{\delta}\mathrm{H}(\delta\rho),	
  \end{multline}
  by substituting $(1-\e)\widetilde{\lambda}^{\min}_{\r}+\epsilon$ for $\lambda$ in \eqref{eq:lmin}. We consolidate notation using $l := 1 - \widetilde{\lambda}^{\min}_{\r}$ and bound the Shannon entropy function from below using the second bound in \eqref{eq:shannon_bounds} to give
  \begin{align}
  & 2\Psi_{\min}\left((1-\e)\widetilde{\lambda}^{\min}_{\r}+\epsilon,\delta,\rho\right) \nonumber\\
  \label{part2r_2b}
  & > 2 \left[-\r\log\r + \r - \r^2\right] + (1-\rho)\log\left[(1-\e)(1-l) + \e\right] + \rho\log\r - \rho  \nonumber\\
  & \quad + 1 - (1-\e)(1-l) - \e + \frac{2}{\d}\left[-\r\log\left(\d\r\right) + \d\r - \d^2\r^2\right],\\
  \label{part2r_3b}
  & = -2\r\log\r + 2\r - 2\r^2 + (1-\rho)\log\left[1 - \e - (1-\e)l + \e \right] + \rho\log\r - \rho \nonumber\\
  & \quad + 1 - 1 + \e + (1-\e)l - \e - 2\r\log\left(\d\r\right) + 2\r - 2\d\r^2,\\
  \label{part2r_4b}
  & = \log\left[1 - (1-\e)l \right] + (1-\e)l - \rho\log\left[1 - (1-\e)l \right] + \rho\log\left(\frac{1}{\d^2\r^3}\right)\nonumber\\
  & \quad + 3\rho - 2(1+\d)\r^2.
  \end{align}
  From \eqref{part2r_2b} to \eqref{part2r_3b} we expanded brackets and simplified and further simplified from \eqref{part2r_3b} to \eqref{part2r_4b}.

  Now replacing $\rho\log\left(\frac{1}{\delta^2\rho^3}\right)$ by $\left(l^2-c\r\right)/2$, bounding above the second logarithmic term using the first bound of \eqref{eq:logm_ubound} and factoring out $\log(1-l)$ we have
  \begin{align}
  & 2\Psi_{\min}\left((1-\e)\widetilde{\lambda}^{\min}_{\r}+\epsilon,\delta,\rho\right) \nonumber\\
  \label{part2r_5b}
  & > \log\left(1 - l \right) + \log\left(1 + \frac{\e l}{1-l} \right) + l - \e l + (1-\e)\rho l + \frac{1}{2}\left(l^2-c\r\right) + 3\rho\nonumber\\
  & \quad - 2(1+\d)\r^2,\\
  \label{part2r_6b}
  & > \log\left(1 - l \right) + l + \frac{1}{2}l^2 - \frac{1}{2}c\r + 3\rho - \e l + \log\left(1 + \e l \right) + \rho l - \e\rho l\nonumber\\
  & \quad - 2(1+\d)\r^2,\\
  \label{part2r_7b}
  & > -l - \frac{1}{2}l^2 - \frac{1}{2}l^3 + l + \frac{1}{2}l^2 + \frac{1}{2}(6-c)\r + \rho l - \e l + \e l - \frac{1}{2}\e^2l^2 - \e\rho l\nonumber\\
  & \quad - 2(1+\d)\r^2,\\
  \label{part2r_8b}
  & = \frac{1}{2}(6-c)\r - \frac{1}{2}l^3 + \rho l - 2(1+\d)\r^2 - \frac{1}{2}\e^2l^2 - \e\rho l.
  \end{align}

  From \eqref{part2r_5b} to \eqref{part2r_6b} we bounded below $1/(1-l)$ using \eqref{eq:r_mfrac_series}. From \eqref{part2r_6b} to \eqref{part2r_7b} we bounded below the first logarithmic term using
  \begin{equation}
  \label{eq:logm_lbound}
  \log(1 - x) \geq -x - \frac{1}{2}x^2 - \frac{1}{2}x^3 \quad \forall x\in[0,0.44],
  \end{equation}
  and also bounded below the second logarithmic term using \eqref{eq:logp_lbound}. From \eqref{part2r_7b} to \eqref{part2r_8b} we simplified.

  The dominant terms in \eqref{part2r_8b} are the first two term, all the rest go to zero faster as $\r\rightarrow 0$. Therefore, for \eqref{part2r_8b} to be positive as $\r\rightarrow 0$ we need the sum of the first two terms to be positive. This means
  \begin{equation}
  \label{eq:fr_p2r_1b}
  \frac{1}{2}(6-c)\r - \frac{1}{2}l^3 > 0 \quad \Rightarrow \quad l^3 < (6-c)\r.
  \end{equation}
  This holds for $c<6$ and small enough $\r$ and since $l$ is a decreasing function of $\r^{-1}$ there would not a $\r_0$ below which this ceases to hold as $\r\rightarrow 0$. Hence we conclude that for fixed $\d, ~\e>0$ and $c<6$ there does not exist a $\r_0$ such that for $\r<\r_0$, \eqref{part2r_8b} is negative and $2\Psi_{\min}\left((1-\e)\widetilde{\lambda}^{\min}_{\r}+\epsilon,\delta,\rho\right) \leq 0$ as $\r\rightarrow 0$.
  Thus
  \begin{multline*}
  Prob\left(L(k,n,N;A) > (1-\e)\widetilde{\mathcal{L}}^{\r}(\delta,\rho)\right) \\ \le poly\left(n,(1-\e)\widetilde{\lambda}^{\min}_{\r}+\epsilon\right)\cdot\exp\left[2n\Psi_{\min}\left((1-\e)\widetilde{\lambda}^{\min}_{\r}+\epsilon,\delta,\rho\right)\right],
  \end{multline*}
   and as $n \rightarrow \infty$ the right hand side of this does not go to zero.

\end{description}

   Now \textbf{Part 1} and \textbf{Part 2} put together shows that $\widetilde{\mathcal{L}}^{\r}(\delta,\rho)$ is also a tight bound of $L(k,n,N;A)$ with overwhelming probability as the problem size grows in the regime prescribed for $\widetilde{\mathcal{L}}^{\r}(\delta,\rho)$ in Theorem \ref{thm:ULbctr}.

\end{proof}

\subsection{Theorem \ref{thm:ULbctd}}\label{sec:pthm:ULbctd}

\subsubsection{The upper bound, $\widetilde{\mathcal{U}}^{\d}(\delta,\rho)$}\label{sec:pUbctdthm}

\begin{proof}
Define $$\widetilde{\lambda}^{\max}_{\d}(\delta,\rho) := 1 + 3\rho + \rho\log\left(\frac{1}{\delta^2\rho^3}\right) + (1+\rho)\log\left[c\log\left(\frac{1}{\delta^2\rho^3}\right)\right].$$ It follows from \eqref{eq:Ubctd} that $\widetilde{\mathcal{U}}^{\d}(\delta,\rho) = \widetilde{\lambda}^{\max}_{\d}(\delta,\rho) - 1$.  Bounding $\widetilde{\mathcal{U}}^{\d}(\delta,\rho)$ above by $\widetilde{\mathcal{U}}^{\d}(\delta,\rho) + \e$ is equivalent to bounding $\widetilde{\lambda}^{\max}_{\d}$ above by $\widetilde{\lambda}^{\max}_{\d}+\epsilon$. We first establish that for a slightly looser bound, with $c>1$, the exponent $\Psi_{\max}\left(\widetilde{\lambda}^{\max}_{\d}+\epsilon,\delta,\rho\right)$ is negative and then verify that when multiplied by $n$ it diverges to $-\infty$ as $n$ increases. We also show that for a slightly tighter bound, with $c\leq \r$, the exponent $\Psi_{\max}\left(\widetilde{\lambda}^{\max}_{\d}-\epsilon,\delta,\rho\right)$ is bounded from below by zero, and hence the bound $\widetilde{\mathcal{U}}^{\d}(\delta,\rho)$ cannot be improved using the inequality \eqref{ppmax} from \cite{blanchard2011compressed}  We show, in two parts that for $\r$ fixed:
\begin{enumerate}
  \item $\exists ~\delta_0, ~\epsilon > 0 ~\rm{and} ~c>1$ such that for $\delta < \delta_0, \Psi_{\max}\left(\widetilde{\lambda}^{\max}_{\d}+\epsilon,\delta,\rho\right) \leq 0;$
  \item $\nexists ~\delta_0, ~\epsilon > 0 ~\rm{and} ~c\leq\r$ such that for $\delta < \delta_0, \Psi_{\max}\left(\widetilde{\lambda}^{\max}_{\d}-\epsilon,\delta,\rho\right) \leq 0.$
\end{enumerate}
which are proven separately in the two parts as follows.

\begin{description}
  \item[Part 1:]
  \begin{multline}
  \label{part1d_1}
  2\Psi_{\max}\left(\widetilde{\lambda}^{\max}_{\d}+\epsilon,\delta,\rho\right) = (1+\rho)\log\left(\widetilde{\lambda}^{\max}_{\d}+\epsilon\right) \\
  - \rho\log(\rho) + \rho + 1 - \left(\widetilde{\lambda}^{\max}_{\d}+\epsilon\right)+ \frac{2}{\delta}\mathrm{H}(\delta\rho),	
  \end{multline}
  by substituting $\widetilde{\lambda}^{\max}_{\r}+\epsilon$ for $\lambda$ in \eqref{eq:lmax}. We bound the Shannon entropy function above using the first bound of \eqref{eq:shannon_bounds} and consolidate notation using $u := \widetilde{\lambda}^{\max}_{\r} - 1$, then \eqref{part1d_1} becomes
  \begin{align}
  &   2\Psi_{\max}\left(\widetilde{\lambda}^{\max}_{\d}+\epsilon,\delta,\rho\right) \qquad\nonumber\\
  \label{part1d_2}
  & < (1+\rho)\log\left[(1+u) + \e\right] - \rho\log\r + \rho + 1 - (1+u) - \e\nonumber\\
  & \quad + \frac{2}{\d}\left[-\d\r\log\left(\d\r\right) + \d\r\right],\\
  \label{part1d_3}
  & = (1+\rho)\log\left(1+u + \e\right) + \r\log\left(\frac{1}{\delta^2\rho^3}\right) + 3\rho - u - \e.
  \end{align}

  From \eqref{part1d_2} to \eqref{part1d_3} we simplified. Next where $u$ is not in the logarithmic term we replace it by $\rho\log\left(\frac{1}{\delta^2\rho^3}\right) + (1+\rho)\log\left[c\log\left(\frac{1}{\delta^2\rho^3}\right)\right] + 3\rho$ to have

  \begin{align}
  &   2\Psi_{\max}\left(\widetilde{\lambda}^{\max}_{\d}+\epsilon,\delta,\rho\right) \nonumber\\
  \label{part1d_4}
  & < (1+\rho)\log\left(1+u + \e\right) + \r\log\left(\frac{1}{\delta^2\rho^3}\right) + 3\rho - \rho\log\left(\frac{1}{\delta^2\rho^3}\right) - 3\rho\nonumber\\
  & \quad - (1+\rho)\log\left[c\log\left(\frac{1}{\delta^2\rho^3}\right)\right] - \e,\\
  \label{part1d_5}
  & = (1+\rho)\log\left(1+u + \e\right) - \e - (1+\rho)\log\left[c\log\left(\frac{1}{\delta^2\rho^3}\right)\right],\\
  \label{part1d_6}
  & = - \a(1+\r) - \e + (1+\rho)\log\left[\frac{1 + u + \e}{c\log\left(\frac{1}{\delta^2\rho^3}\right)} \right] + \a(1+\r),\qquad\\
  \label{part1d_7}
  & = -\a - \a\r - \e + (1+\rho)\log\left[\frac{1 + u + \e}{c\log\left(\frac{1}{\delta^2\rho^3}\right)} \right] + \a(1+\r)\log e,\qquad\\
  \label{part1d_8}
  & < -\a + (1+\rho)\log\left[\frac{e^{\a}(1 + u + \e)}{c\log\left(\frac{1}{\delta^2\rho^3}\right)} \right].
  \end{align}

  From \eqref{part1d_4} to \eqref{part1d_5} we simplified and from \eqref{part1d_5} to \eqref{part1d_6} we combined the logarithmic terms and to create a constant we add $-\a(1+\r)$ and $\a(1+\r)$ for a small positive constant $0<\a<1$. From \eqref{part1d_6} to \eqref{part1d_7} we rewrote $\a(1+\r)$ as $\a(1+\r)\log e$. From \eqref{part1d_7} to \eqref{part1d_8} incorporated the second logarithmic term into the first one and we bounded above \eqref{part1d_7} by dropping the $-\e$ and $-\a\r$.

  Equation \eqref{part1d_8} is clearly negative if the second term is negative, which is satisfied if the argument of the logarithm to be less than one. This leads to
  \begin{equation}
  \label{eq:fr_p1d_1}
  e^{-\a}c\log\left(\frac{1}{\delta^2\rho^3}\right) \geq u + 1 + \e,
  \end{equation}
  where again substituting $\rho\log\left(\frac{1}{\delta^2\rho^3}\right) + (1+\rho)\log\log\left(\frac{1}{\delta^2\rho^3}\right) + 3\rho$ for $u$ and reordering the right hand side of \eqref{eq:fr_p1d_1} gives
  \begin{multline}
  \label{eq:fr_p1d_2}
   e^{-\a}c\log\left(\frac{1}{\delta^2\rho^3}\right) \geq \log\log\left(\frac{1}{\delta^2\rho^3}\right)  + 1 + \e \\ + \r \left[3 + \log\left(\frac{1}{\delta^2\rho^3}\right) + \log\log\left(\frac{1}{\delta^2\rho^3}\right)\right].
  \end{multline}
  For small $0<\a<1$ and $c>1$, the left hand side of \eqref{eq:fr_p1d_2} is an unbounded strictly increasing function of $\d^{-1}$ growing exponentially faster than the right hand side of \eqref{eq:fr_p1d_2}.  Consequently there is a unique $\d_0$ for which the inequality \eqref{eq:fr_p1d_2} holds for fixed $\r, ~\e>0, ~c>1$ and any $\d\le\d_0$ and as a result making $2\Psi_{\max}\left(\widetilde{\lambda}^{\max}_{\d}+\epsilon,\delta,\rho\right) < 0$.

  Having established a negative bound from above and the $\d_0$ for which it is valid, it remains to show that $~n \cdot 2\Psi_{\max}\left(\widetilde{\lambda}^{\max}_{\d}+\epsilon,\delta,\rho\right) \rightarrow -\infty$ as $(k,n,N) \rightarrow \infty$, which verifies an exponential decay to zero of the bound \eqref{ppmax} with $n$. This follows from the first term of the right hand side of \eqref{part1d_8}, giving a concluding bound $n(-\a).$ For $\d<\d_0$ therefore
  \begin{multline*}
  Prob\left(U(k,n,N;A) > \widetilde{\mathcal{U}}^{\d}(\delta,\rho) + \e\right) \\ \le poly\left(n,\widetilde{\lambda}^{\max}_{\d}+\epsilon\right)\cdot\exp\left(-\a n\right).
  \end{multline*}
  The right hand side of which goes to zero as $n\rightarrow \infty$.\\

  \item[Part 2:]
  \begin{multline}
  \label{part2d_1}
  2\Psi_{\max}\left(\widetilde{\lambda}^{\max}_{\d}-\epsilon,\delta,\rho\right) = (1+\rho)\log\left(\widetilde{\lambda}^{\max}_{\d}-\epsilon\right) \\
  - \rho\log(\rho) + \rho + 1 - \left(\widetilde{\lambda}^{\max}_{\d}-\epsilon\right)+ \frac{2}{\delta}\mathrm{H}(\delta\rho),	
  \end{multline}
  by substituting $\widetilde{\lambda}^{\max}_{\r}-\epsilon$ for $\lambda$ in \eqref{eq:lmax}. We lower bound the Shannon entropy function using the second bound of \eqref{eq:shannon_bounds} and consolidate notation using $u := \widetilde{\lambda}^{\max}_{\d} - 1$, then \eqref{part2d_1} becomes
  \begin{align}
  &   2\Psi_{\max}\left(\widetilde{\lambda}^{\max}_{\d}-\epsilon,\delta,\rho\right) \nonumber\\
  \label{part2d_2}
  & > (1+\rho)\log\left[(1+u) - \e\right] - \rho\log\r + \rho + 1 - (1+u) + \e\nonumber\\
  & \quad - 2\r\log\left(\d\r\right) + \frac{2}{\d}\left[-\d\r\log\left(\d\r\right) + \d\r - \d^2\r^2\right],\\
  \label{part2d_3}
  & = (1+\rho)\log\left(u+1-\e\right) + \r\log\left(\frac{1}{\delta^2\rho^3}\right) + 3\rho - u + \e - 2\d\r^2,\\
  \label{part2d_4}
  & = (1+\rho)\log\left(u+1-\e\right) + \r\log\left(\frac{1}{\delta^2\rho^3}\right) + 3\rho - \rho\log\left(\frac{1}{\delta^2\rho^3}\right) - 3\rho\nonumber\\
  & \quad - (1+\rho)\log\left[c\log\left(\frac{1}{\delta^2\rho^3}\right)\right] + \e - 2\d\r^2,\\
  \label{part2d_5}
  & = (1+\rho)\log\left(u+1-\e\right) - (1+\rho)\log\left[c\log\left(\frac{1}{\delta^2\rho^3}\right)\right] + \e - 2\d\r^2,\\
  \label{part2d_6}
  & = \e + (1+\rho)\log\left[\frac{1 + u - \e}{c\log\left(\frac{1}{\delta^2\rho^3}\right)} \right] - 2\d\r^2.
  \end{align}

  From \eqref{part2d_2} to \eqref{part2d_3} we simplified. Then from \eqref{part2d_3} to \eqref{part2d_4} we replace $u$ by $\rho\log\left(\frac{1}{\delta^2\rho^3}\right) + (1+\rho)\log\left[c\log\left(\frac{1}{\delta^2\rho^3}\right)\right] + 3\rho$ where $u$ is not in the logarithmic term. From \eqref{part2d_4} to \eqref{part2d_5} we simplified and from \eqref{part2d_5} to \eqref{part2d_6} we combined the logarithmic terms.

  The last term in \eqref{part2d_6} obviously goes to zero as $\d\rightarrow 0$, then for the expression to remain positive we need to know how the dominant term, which is the second term, behaves. For this term to be nonnegative as $\d\rightarrow 0$ for fixed $\r$ we need the argument of the logarithmic to be greater than or equal to 1 which means the following.
  \begin{equation*}
  u + 1 + \e \geq c\log\left(\frac{1}{\delta^2\rho^3}\right).
  \end{equation*}
  Therefore substituting for $u$ we have
  \begin{equation*}
  \r\log\left(\frac{1}{\delta^2\rho^3}\right) + (1+\r)\log\left[c\log\left(\frac{1}{\delta^2\rho^3}\right)\right] + 3\r  + 1 + \e \geq c\log\left(\frac{1}{\delta^2\rho^3}\right),
  \end{equation*}
 Then we expand the second logarithmic term and rearrange to get
  \begin{equation}
   \label{eq:fr_p2d_1}
  (\r-c)\log\left(\frac{1}{\delta^2\rho^3}\right) + (1+\r)\log\left[c\log\left(\frac{1}{\delta^2\rho^3}\right)\right] + 3\r  + 1 + \e \geq 0.
  \end{equation}

  Inequality \eqref{eq:fr_p2d_1} is always true for fixed $\r$ and $c<\r$ as $\d \rightarrow 0$. Therefore, we conclude that there does not exists $\d_0$ such that for any $\r$ fixed and $\e>0$ for $\d<\d_0$ \eqref{part2d_6} is negative and $2\Psi_{\max}\left(\widetilde{\lambda}^{\max}_{\d}-\epsilon,\delta,\rho\right) < 0$ as $\d\rightarrow 0$. Thus
  \begin{multline*}
  Prob\left(U(k,n,N;A) > \widetilde{\mathcal{U}}^{\d}(\delta,\rho)-\e\right) \\ \le poly\left(n,\widetilde{\lambda}^{\max}_{\d}-\epsilon\right)\cdot\exp\left[2n\Psi_{\max}\left(\widetilde{\lambda}^{\max}_{\d}-\epsilon,\delta,\rho\right)\right],
  \end{multline*}
   and as $n \rightarrow \infty$ the right hand side of this does not necessarily go to zero.

\end{description}

Now \textbf{Part 1} and \textbf{Part 2} put together shows that $\widetilde{\mathcal{U}}^{\d}(\delta,\rho)$ is also a tight upper bound of $U(k,n,N;A)$ with overwhelming probability as the problem size grows in the regime prescribed for $\widetilde{\mathcal{U}}^{\d}(\delta,\rho)$ in Theorem \ref{thm:ULbctd}.

\end{proof}

\subsubsection{The lower bound, $\widetilde{\mathcal{L}}^{\d}(\delta,\rho)$}\label{sec:pLbctdthm}

\begin{proof}
Define $$\widetilde{\lambda}^{\min}_{\d}(\delta,\rho) := \exp\left(-\frac{3\rho+c}{1-\rho}\right) \cdot \left(\delta^2\rho^3\right)^{\frac{\rho}{1-\rho}}, \quad \Rightarrow \quad \widetilde{\mathcal{L}}^{\d}(\delta,\rho) = 1 - \widetilde{\lambda}^{\min}_{\d}(\delta,\rho)$$ as from \eqref{eq:Lbctd}. Bounding $\widetilde{\mathcal{L}}^{\d}(\delta,\rho)$ above by $(1+\e)\widetilde{\mathcal{L}}^{\d}(\delta,\rho)$ is equivalent to bounding $\widetilde{\lambda}^{\min}_{\d}$ above by $(1+\e)\widetilde{\lambda}^{\min}_{\d}-\epsilon$. We first establish for a slightly looser bound, with $c>1$, the exponent  $\Psi_{\min}\left((1+\e)\widetilde{\lambda}^{\min}_{\d}-\epsilon,\delta,\rho\right)$ is negative and then verify that when multiplied by $n$ it diverges to $-\infty$ as $n$ increases. We also show that for a slightly tighter bound, with $c<1$, the exponent $\Psi_{\min}\left((1-\e)\widetilde{\lambda}^{\min}_{\d}+\epsilon,\delta,\rho\right)$ is bounded from below by zero, and hence the bound $\widetilde{\mathcal{L}}^{\d}(\delta,\rho)$ cannot be improved using the inequality \eqref{ppmin} from \cite{blanchard2011compressed}. We show, in two parts that for $\r$ fixed:
\begin{enumerate}
  \item $\exists ~\d_0, ~\epsilon > 0 ~\rm{and} ~c>1$ such that for $\d < \d_0, \Psi_{\min}\left((1+\e)\widetilde{\lambda}^{\min}_{\d}-\epsilon,\delta,\rho\right) \leq 0;$
  \item $\nexists ~\d_0, ~\epsilon > 0 ~\rm{and} ~c<1$ such that for $\d < \d_0, \Psi_{\min}\left((1-\e)\widetilde{\lambda}^{\min}_{\d}+\epsilon,\delta,\rho\right) \leq 0,$
\end{enumerate}
which are proven separately in the two parts as follows.

\begin{description}
  \item[Part 1:]
  \begin{multline}
  \label{part1d_1b}
  2\Psi_{\min}\left((1+\e)\widetilde{\lambda}^{\min}_{\d}-\epsilon,\delta,\rho\right) = 2\mathrm{H}(\rho) + (1-\rho)\log\left((1+\e)\widetilde{\lambda}^{\min}_{\d}-\epsilon\right)\\ + \rho\log(\rho) - \rho + 1 - \left((1+\e)\widetilde{\lambda}^{\min}_{\d}-\epsilon\right)+ \frac{2}{\delta}\mathrm{H}(\delta\rho),	
  \end{multline}
  by substituting $(1+\e)\widetilde{\lambda}^{\min}_{\d}-\epsilon$ for $\lambda$ in \eqref{eq:lmin}. We now upper bound the Shannon entropy terms using the first bound of \eqref{eq:shannon_bounds} and factor out $\widetilde{\lambda}^{\min}_{\d}$ for \eqref{part1d_1b} to become
  \begin{align}
  & 2\Psi_{\min}\left((1+\e)\widetilde{\lambda}^{\min}_{\d}-\epsilon,\delta,\rho\right)\nonumber\\
  \label{part1d_2b}
  & < 2\left[-\r\log\r + \r - \r^2\right] + (1-\rho)\log\left(\widetilde{\lambda}^{\min}_{\d}\right)  - (1+\e)\widetilde{\lambda}^{\min}_{\d} + \epsilon + 1 - \rho \nonumber\\
  & \quad + \rho\log\r + (1-\rho)\log\left[\frac{(1+\e)\widetilde{\lambda}^{\min}_{\d}-\epsilon}{\widetilde{\lambda}^{\min}_{\d}}\right] + \frac{2}{\d}\left[-\r\log\left(\d\r\right) + \d\r\right],\\
  \label{part1d_3b}
  & = (1-\rho)\log\left(\widetilde{\lambda}^{\min}_{\d}\right) - (1+\e)\widetilde{\lambda}^{\min}_{\d} + \epsilon + (1-\rho)\log\left[(1+\e)-\frac{\epsilon}{\widetilde{\lambda}^{\min}_{\d}}\right]\nonumber\\
  & \quad + \rho\log\left(\frac{1}{\delta^2\rho^3}\right) + 3\rho + 1.
  \end{align}
  From \eqref{part1d_2b} to \eqref{part1d_3b} we simplified. Using the fact that by the definition of $\widetilde{\mathcal{L}}^{\d}(\delta,\rho)$ in \eqref{eq:Lbctd} $$\log\left(\widetilde{\lambda}^{\min}_{\d}\right) = -\frac{\rho}{1-\rho}\log\left(\frac{1}{\delta^2\rho^3}\right) - \frac{3\rho + c}{1-\rho},$$ we substitute this in \eqref{part1d_3b} for $\log\left(\widetilde{\lambda}^{\min}_{\d}\right)$ to get
  \begin{align}
  & 2\Psi_{\min}\left((1+\e)\widetilde{\lambda}^{\min}_{\d}-\epsilon,\delta,\rho\right)\nonumber\\
  \label{part1d_4b}
  & < (1-\rho)\left[-\frac{\rho}{1-\rho}\log\left(\frac{1}{\delta^2\rho^3}\right) - \frac{3\rho + c}{1-\rho}\right] - (1+\e)\widetilde{\lambda}^{\min}_{\d} + \epsilon\nonumber\\
  & \quad + (1-\rho)\log\left[(1+\e)-\frac{\epsilon}{\widetilde{\lambda}^{\min}_{\d}}\right] + \rho\log\left(\frac{1}{\delta^2\rho^3}\right) + 3\rho + 1,\\
  \label{part1d_5b}
  & = -\r\log\left(\frac{1}{\delta^2\rho^3}\right) - 3\rho - c - (1+\e)\widetilde{\lambda}^{\min}_{\d} + \epsilon + \rho\log\left(\frac{1}{\delta^2\rho^3}\right)\nonumber\\
  & \quad + (1-\rho)\log\left[(1+\e)-\frac{\epsilon}{\widetilde{\lambda}^{\min}_{\d}}\right] + 3\rho + 1,\\
  \label{part1d_6b}
  & = (1-\rho)\log\left[(1+\e) - \frac{\epsilon}{\widetilde{\lambda}^{\min}_{\d}}\right] - \widetilde{\lambda}^{\min}_{\d} - \e\widetilde{\lambda}^{\min}_{\d} - (c-1) + \epsilon.
  \end{align}
  From \eqref{part1d_4b} to \eqref{part1d_5b} we expanded the brackets and from \eqref{part1d_5b} to \eqref{part1d_6b} we simplified. Now we consolidate notation using $l := 1 - \widetilde{\lambda}^{\min}_{\d}$ and substituting this in \eqref{part1d_6b} we have
  \begin{align}
  &2\Psi_{\min}\left((1+\e)\widetilde{\lambda}^{\min}_{\d}-\epsilon,\delta,\rho\right) \nonumber\\
  \label{part1d_7b}
  & < (1-\rho)\log\left[(1+\e) - \frac{\epsilon}{1 - l}\right] - (1 - l) - \e(1-l)  - (c-1) + \epsilon,\\
  \label{part1d_8b}
  & = - (c-1) + (1-\rho)\log\left(1 - \frac{\epsilon l}{1 - l}\right) - (1 - l) + \e l,\\
  \label{part1d_9b}
  & < - (c-1) + \e l  - (1-\rho)\frac{\epsilon l}{1 - l} - (1-l),\\
  \label{part1d_10b}
  & = - \frac{1}{2}(c-1) - \frac{1}{2}(c-1) + \e l.
  \end{align}

  From \eqref{part1d_7b} to \eqref{part1d_8b} we simplified and from \eqref{part1d_8b} to \eqref{part1d_9b} we bounded above the logarithmic term using the first bound of \eqref{eq:logm_ubound}. From \eqref{part1d_9b} to \eqref{part1d_10b} we dropped the third and fourth terms, which are negative, and split the leading term into half. Inequality \eqref{part1d_10b} can be further bounded by $-(c-1)/2$ (which will be negative if $c>1$) by choosing $\e$ to be less than $(c-1)/2$ and noting that $l\in (0,1]$.

  Having established a negative bound from above and the $\d_0$ for which it is valid, it remains to show that $~n \cdot 2\Psi_{\min}\left((1+\e)\widetilde{\lambda}^{\min}_{\d}-\epsilon,\delta,\rho\right) \rightarrow -\infty$ as $(k,n,N) \rightarrow \infty$, which verifies an exponential decay to zero of the bound \eqref{ppmin} with $n$. This follows from the first term of the right hand side of \eqref{part1d_10b} giving a concluding bound $-\frac{1}{2}(c-1)n.$ For $\d<\d_0$ therefore
  \begin{multline*}
  Prob\left(L(k,n,N;A) > (1+\e)\widetilde{\mathcal{L}}^{\d}(\delta,\rho)\right) \\ \le poly\left(n,(1+\e)\widetilde{\lambda}^{\min}_{\d}-\epsilon\right)\cdot\exp\left[-\frac{(c-1)n}{2}\right].
  \end{multline*}
  The right hand side of which goes to zero as $n\rightarrow \infty$.\\

  \item[Part 2:]
  \begin{multline}
  \label{part2d_1b}
  2\Psi_{\min}\left((1-\e)\widetilde{\lambda}^{\min}_{\d}+\epsilon,\delta,\rho\right) = 2\mathrm{H}(\rho) + (1-\rho)\log\left((1-\e)\widetilde{\lambda}^{\min}_{\d}+\epsilon\right)\\ + \rho\log(\rho) - \rho + 1 - \left((1-\e)\widetilde{\lambda}^{\min}_{\d}+\epsilon\right)+ \frac{2}{\delta}\mathrm{H}(\delta\rho),	
  \end{multline}
  by substituting $(1-\e)\widetilde{\lambda}^{\min}_{\d}+\epsilon$ for $\lambda$ in \eqref{eq:lmin}. Next we bound the Shannon entropy functions from below using the second bound in \eqref{eq:shannon_bounds} to give

  \begin{align}
  & 2\Psi_{\min}\left((1-\e)\widetilde{\lambda}^{\min}_{\d}+\epsilon,\delta,\rho\right) \nonumber\\
  \label{part2d_2b}
  & > 2\left[-\r\log\r + \r - \r^2\right] + (1-\rho)\log\left(\widetilde{\lambda}^{\min}_{\d}\right)  - (1-\e)\widetilde{\lambda}^{\min}_{\d} + 1 - \epsilon \nonumber\\
  & \quad + \rho\log\r - \rho + (1-\rho)\log\left[\frac{(1-\e)\widetilde{\lambda}^{\min}_{\d}+\epsilon}{\widetilde{\lambda}^{\min}_{\d}}\right]\nonumber\\
  & \quad + \frac{2}{\d}\left[-\r\log\left(\d\r\right) + \d\r - \d^2\r^2\right],\\
  \label{part2d_3b}
  & = (1-\rho)\log\left(\widetilde{\lambda}^{\min}_{\d}\right) - (1-\e)\widetilde{\lambda}^{\min}_{\d} - \epsilon + (1-\rho)\log\left[(1-\e)+\frac{\epsilon}{\widetilde{\lambda}^{\min}_{\d}}\right]\nonumber\\
  & \quad + \rho\log\left(\frac{1}{\delta^2\rho^3}\right) + 3\rho + 1 - 2(1+\d)\r^2.
  \end{align}

   From \eqref{part2d_2b} to \eqref{part2d_3b} we simplified. Using the fact that by the definition of $\widetilde{\mathcal{L}}^{\d}(\delta,\rho)$ in \eqref{eq:Lbctd} $$\log\left(\widetilde{\lambda}^{\min}_{\d}\right) = -\frac{\rho}{1-\rho}\log\left(\frac{1}{\delta^2\rho^3}\right) - \frac{3\rho+c}{1-\rho},$$ we substitute this in \eqref{part2d_3b} for $\log\left(\widetilde{\lambda}^{\min}_{\d}\right)$ to get

  \begin{align}
  & 2\Psi_{\min}\left((1-\e)\widetilde{\lambda}^{\min}_{\d}+\epsilon,\delta,\rho\right)\nonumber\\
  \label{part2d_4b}
  & > (1-\rho)\left[-\frac{\rho}{1-\rho}\log\left(\frac{1}{\delta^2\rho^3}\right) - \frac{3\rho+c}{1-\rho}\right] - (1-\e)\widetilde{\lambda}^{\min}_{\d} - \epsilon + 3\rho\nonumber\\
  & \quad + (1-\rho)\log\left[(1-\e)+\frac{\epsilon}{\widetilde{\lambda}^{\min}_{\d}}\right] + \rho\log\left(\frac{1}{\delta^2\rho^3}\right) + 1 - 2(1+\d)\r^2,\\
  \label{part2d_5b}
  & = -\r\log\left(\frac{1}{\delta^2\rho^3}\right) - 3\rho - c - (1-\e)\widetilde{\lambda}^{\min}_{\d} - \epsilon + \rho\log\left(\frac{1}{\delta^2\rho^3}\right)\nonumber\\
  & \quad + (1-\rho)\log\left[(1-\e)+\frac{\epsilon}{\widetilde{\lambda}^{\min}_{\d}}\right] + 3\rho + 1 - 2(1+\d)\r^2,\\
  \label{part2d_6b}
  & = (1-\rho)\log\left[(1-\e) + \frac{\epsilon}{\widetilde{\lambda}^{\min}_{\d}}\right] - \widetilde{\lambda}^{\min}_{\d} + \e\widetilde{\lambda}^{\min}_{\d} - \epsilon + 1 - c\nonumber\\
  & \quad - 2(1+\d)\r^2.
  \end{align}

  From \eqref{part2d_4b} to \eqref{part2d_5b} we expanded the brackets and from \eqref{part2d_5b} to \eqref{part2d_6b} we simplified. Now we consolidate notation using $l := 1 - \widetilde{\lambda}^{\min}_{\d}$ and substituting this in \eqref{part2d_6b} we have
  \begin{align}
  &2\Psi_{\min}\left((1-\e)\widetilde{\lambda}^{\min}_{\d}+\epsilon,\delta,\rho\right) \nonumber\\
  \label{part2d_7b}
  & > (1-\rho)\log\left[(1-\e) + \frac{\e}{1-l}\right] - (1 - l) + \e(1 - l) - \epsilon + 1 - c\nonumber\\
  & \quad - 2(1+\d)\r^2,\\
  \label{part2d_8b}
  & = (1-\rho)\log\left(1 + \frac{\e l}{1 - l}\right) + l - c - \e l - 2(1+\d)\r^2,\\
  \label{part2d_9b}
  & > (1-\rho)\log\left(1 + \e l\right) + l - c -  \e l - 2(1+\d)\r^2,\\
  \label{part2d_10b}
  & > (1-\rho)\left(\e l - \frac{1}{2} \e^2l^2\right) + l - c -  \e l - 2(1+\d)\r^2,\\
  \label{part2d_11b}
  & = \e l - \frac{1}{2} \e^2l^2  - \e\rho l + \frac{1}{2} \e^2\r l^2 + l - c -  \e l - 2\r^2 - 2\d\r^2,\\
  \label{part2d_12b}
  & = l - c  - 2\r^2 - \e l - \e\rho l - \frac{1}{2} \e^2l^2  + \frac{1}{2} \e^2\r l^2 - 2\d\r^2.
  \end{align}
  We simplified from \eqref{part2d_7b} to \eqref{part2d_8b} and from \eqref{part2d_8b} to \eqref{part2d_9b} we bounded below $1/(1-l)$ using the bound of \eqref{eq:r_mfrac_series}. From \eqref{part2d_9b} to \eqref{part2d_10b} we bounded below the logarithmic term using the bound of \eqref{eq:logp_lbound}. From \eqref{part2d_10b} to \eqref{part2d_11b} we expanded the brackets and from \eqref{part2d_11b} to \eqref{part2d_12b} we simplified.

  The leading terms of \eqref{part2d_12b} are the first three and $l$ is strictly increasing as $\d^{-1}$ approaches 1. If $c<1$, there will be some values of $\r$ for which \eqref{part2d_12b} will always be positive as $\d\rightarrow 0$. Thus there does not exist any $\d_0$ such that for any $\r$ fixed, $\e>0, ~c<1$ and $\d<\d_0$, \eqref{part2d_12b} becomes negative. Thus
  \begin{multline*}
  Prob\left(L(k,n,N;A) > (1-\e)\widetilde{\mathcal{L}}^{\d}(\delta,\rho)\right) \\ \le poly\left(n,(1-\e)\widetilde{\lambda}^{\min}_{\d}+\epsilon\right)\cdot\exp\left[2n\Psi_{\min}\left((1-\e)\widetilde{\lambda}^{\min}_{\d}+\epsilon,\delta,\rho\right)\right],
  \end{multline*}
   and as $n \rightarrow \infty$ the right hand side of this does not necessarily go to zero.

\end{description}
Now \textbf{Part 1} and \textbf{Part 2} put together shows that $\widetilde{\mathcal{L}}^{\d}(\delta,\rho)$ is also a tight bound of $L(k,n,N;A)$ with overwhelming probability as the sample size grows in the regime prescribed for $\widetilde{\mathcal{L}}^{\d}(\delta,\rho)$ in Theorem \ref{thm:ULbctd}.

\end{proof}

\subsection{Theorem \ref{thm:ULbctg}}\label{sec:pthm:ULbctg}

\subsubsection{The upper bound, $\widetilde{\mathcal{U}}^{\g}(\delta,\rho_{\g}(\d))$}\label{sec:pUbctgthm}

\begin{proof}
To simplify notation we will use $\r$ for $\rho_{\g}(\d)$ in the proof. Lets define $$\widetilde{\lambda}^{\max}_{\g}(\delta,\rho) := 1 + \sqrt{2\rho\log\left(\frac{1}{\delta^2\rho^3}\right) + 6\r} + c_u\left[2\rho\log\left(\frac{1}{\delta^2\rho^3}\right) + 6\r\right].$$ It follows from \eqref{eq:Ubctg} that $\widetilde{\mathcal{U}}^{\g}(\delta,\rho) = \widetilde{\lambda}^{\max}_{\g}(\delta,\rho) - 1$. Bounding $\widetilde{\mathcal{U}}^{\g}(\delta,\rho)$ above by $\widetilde{\mathcal{U}}^{\g}(\delta,\rho)+\e$ is equivalent to bounding $\widetilde{\lambda}^{\max}_{\g}$ above by $\widetilde{\lambda}^{\max}_{\g}+\epsilon$. We first establish that for a slightly looser bound, with $c_u>1/3$, the exponent $\Psi_{\max}\left(\widetilde{\lambda}^{\max}_{\g}+\epsilon,\delta,\rho\right)$ is negative and then verify that when multiplied by $n$ it diverges to $-\infty$ as $n$ increases. We also show that for a slightly tighter bound, with $c_u\leq 1/5$, the exponent $\Psi_{\max}\left(\widetilde{\lambda}^{\max}_{\g}-\epsilon,\delta,\rho\right)$ is bounded from below by zero, and hence the bound $\widetilde{\mathcal{U}}^{\g}(\delta,\rho)$ cannot be improved using the inequality \eqref{ppmax} from \cite{blanchard2011compressed}. We show, in two parts that for $\g > \g_0$ fixed:
\begin{enumerate}
  \item $\exists ~\delta_0, ~\epsilon > 0 ~\rm{and} ~c_u>1/3$ such that for $\delta < \delta_0, \Psi_{\max}\left(\widetilde{\lambda}^{\max}_{\g}+\epsilon,\delta,\rho\right) \leq 0;$
  \item $\nexists ~\delta_0, ~\epsilon > 0 ~\rm{and} ~c_u\leq1/5$ such that for $\delta < \delta_0, \Psi_{\max}\left(\widetilde{\lambda}^{\max}_{\g}-\epsilon,\delta,\rho\right) \leq 0.$
\end{enumerate}
which are proven separately in the two parts.

\begin{description}
  \item[Part 1:]
  \begin{multline}
  \label{part1g_1}
  2\Psi_{\max}\left(\widetilde{\lambda}^{\max}_{\g}+\epsilon,\delta,\rho\right) = \left(1+\r\right)\log\left(\widetilde{\lambda}^{\max}_{\g}+\epsilon\right) - \r\log\left(\r\right)\\ + \r + 1 - \widetilde{\lambda}^{\max}_{\g} - \epsilon + \frac{2}{\delta}H\left(\delta\r\right),
  \end{multline}
  by substituting $\widetilde{\lambda}^{\max}_{\g}+\epsilon$ for $\lambda$ in the definition of $\Psi_{\max}\left(\lambda,\delta,\rho\right)$ in \eqref{eq:lmax}.

  Now letting $u = \widetilde{\lambda}^{\max}_{\g} - 1$ and substituting this in \eqref{part1g_1} and upper bounding the Shannon entropy term using the first bound of \eqref{eq:shannon_bounds} gives \eqref{part1g_2} below
  \begin{align}
  & 2\Psi_{\max}\left(\widetilde{\lambda}^{\max}_{\g}+\epsilon,\delta,\rho\right) \nonumber\\
  \label{part1g_2}
  & < \left(1+\r\right)\log\left(1+u+\epsilon\right) - \r\log\left(\r\right) + \r + 1 - (1+u) - \epsilon \qquad \nonumber\\
  & \quad + \frac{2}{\d}\left[-\d\r\log\left(\d\r\right) + \d\r\right],\\
  \label{part1g_3}
  & = \log(1 + u + \e) + \r\log(1 + u+\e) - u - \epsilon + \r\log\left(\frac{1}{\d^2\r^3}\right) + 3\r,\\
  \label{part1g_4}
  & = \log(1 + u) + \log\left(1 + \frac{\e}{1+u}\right) + \r\log(1 + u+\e) - u - \epsilon \qquad\nonumber\\
  & \quad + \r\log\left(\frac{1}{\d^2\r^3}\right) + 3\r,\\
  \label{part1g_5}
  & < - u + u - \frac{1}{2}u^2 + \frac{1}{3}u^3 + \r\log\left(\frac{1}{\d^2\r^3}\right) + 3\r + \r\log(1 + u+\e) - \epsilon \qquad\nonumber\\
  & \quad + \log\left(1 + \e\right),\\
  \label{part1g_6}
  & < - \frac{1}{2}u^2 + \frac{1}{3}u^3 + \r\log\left(\frac{1}{\d^2\r^3}\right) + 3\r + \r\log(1 + u+\e) - \epsilon + \e.
  \end{align}
  From \eqref{part1g_2} to \eqref{part1g_3} we expanded the $\left(1+\r\right)$ in the first term and simplified while from \eqref{part1g_3} to \eqref{part1g_4} we expanded the first logarithmic term. From \eqref{part1g_4} to \eqref{part1g_5} we bounded above $\log(1+u)$ and $1/(1+u)$ using the second bound of \eqref{eq:logp_ubound} and the bound of \eqref{eq:r_pfrac_series} respectively. Then from \eqref{part1g_5} to \eqref{part1g_6} we simplified and bounded above $\log(1+\e)$ using the first bound of \eqref{eq:logp_ubound}.

  Let $x = 2\r\log\left(\frac{1}{\d^2\r^3}\right) + 6\r$ which means $u = \sqrt{x} + c_ux$. We simplify \eqref{part1g_6} and replace the sum of the second two terms by $\frac{1}{2}x$ and $u$ in the first two terms by $\sqrt{x} + c_ux$ to get
  \begin{align}
  & 2\Psi_{\max}\left(\widetilde{\lambda}^{\max}_{\g}+\epsilon,\delta,\rho\right) \nonumber\\
  \label{part1g_7}
  & < - \frac{1}{2}\left(\sqrt{x} + c_ux\right)^2 + \frac{1}{3}\left(\sqrt{x} + c_ux\right)^3 + \frac{1}{2}x + \r\log(1 + u+\e),\\
  \label{part1g_8}
  & = -\frac{1}{2}x - c_ux^{3/2} - \frac{1}{2}c_u^2x^2 + \frac{1}{3}x^{3/2} + c_ux^2 + c_u^2x^{5/2} + \frac{1}{3}c_u^3x^3 + \frac{1}{2}x \qquad\nonumber\\
  & \quad + \r\log(1 + u+\e),\\
  \label{part1g_9}
  & = - \left(c_u - \frac{1}{3}\right)x^{3/2}  + c_ux^2 - \frac{1}{2}c_u^2x^2 + c_u^2x^{5/2} + \frac{1}{3}c_u^3x^3  \qquad\nonumber\\
  & \quad + \r\log(1 + u+\e).
  \end{align}

  From \eqref{part1g_7} to \eqref{part1g_8} we expanded the first two brackets and from \eqref{part1g_8} to \eqref{part1g_9} we simplified. Substituting $1/\left[\g\log\left(\frac{1}{\d}\right)\right]$ for $\r$ in the expression for $x$ we have $x = 4/\g + g(\r)$ where $g(\r) = 6\r\log\left(1/\r\right) + 6\r$ and goes to zero with $\d$. Therefore, if $4/\g<1$ for $\d$ small enough we will have $x<1$. This means for $\g > 4$ we can define $\d_1$ such that for $\d<\d_1, ~x<1$ and we can upper bound $x^{5/2}$ and $x^3$ by $x^2$ since $x^2 > x^{2+j}$ for $j>0$ when $x<1$. Using this fact we can bound \eqref{part1g_9} above to get
  \begin{align}
  & 2\Psi_{\max}\left(\widetilde{\lambda}^{\max}_{\g}+\epsilon,\delta,\rho\right) \nonumber\\
  \label{part1g_10}
  & < - \left(c_u - \frac{1}{3}\right)x^{3/2} + c_ux^2 - \frac{1}{2}c_u^2x^2 + c_u^2x^2 + \frac{1}{3}c_u^3x^2  \qquad \quad\nonumber\\
  & \quad + \r\log(1 + u+\e),\\
  \label{part1g_11}
  & = - \frac{1}{2}\left(c_u - \frac{1}{3}\right)x^{3/2}  - \frac{1}{2}\left(c_u - \frac{1}{3}\right)x^{3/2} + c_ux^2 + \frac{1}{2}c_u^2x^2 + \frac{1}{3}c_u^3x^2\nonumber\\
  & \quad + \r\log(1 + u+\e).
  \end{align}
  From \eqref{part1g_10} to \eqref{part1g_11} we simplified and split the first term into half. The last term goes to zero with $\d$ so we can define $\d_2$ such that for $\d<\d_2$ we can bound this term above by $x^2$. But also $x^{3/2} = 8/\sqrt{\g^3} + G(\r)$ where $G(\r)$ is the difference between $\left[4/\g + g(\r)\right]^{3/2}$ and $\left(4/\g\right)^{3/2}$ which also goes to zero with $\d$ because this difference is a sum of products with $g(\r)$. This means $-x^{3/2} < -8/\sqrt{\g^3}$ since $g(\r)$ is positive. Now let $f_u(c_u) = c_u + c_u^2/2 + c_u^3/3$, which is positive for all $c_u>0$, using the above therefore we can bound \eqref{part1g_11} to get
  \begin{align}
  & 2\Psi_{\max}\left(\widetilde{\lambda}^{\max}_{\g}+\epsilon,\delta,\rho\right) \nonumber\\
  \label{part1g_12}
  & < \frac{1}{2}\left(c_u - \frac{1}{3}\right)\cdot \left(-\frac{8}{\sqrt{\g^3}}\right)  - \frac{1}{2}\left(c_u - \frac{1}{3}\right)x^{3/2} + f_u(c_u)x^2 + x^2,\\
  \label{part1g_13}
  & = -\frac{4}{\sqrt{\g^3}}\left(c_u - \frac{1}{3}\right)  - \frac{1}{2}\left(c_u - \frac{1}{3}\right)x^{3/2} + \left[1 + f_u(c_u)\right]x^2.
  \end{align}

  From \eqref{part1g_12} to \eqref{part1g_13} we simplified. For \eqref{part1g_13} to be negative all we need is for $c_u>1/3$ and the sum of the last two terms to be non positive, that is:
  \begin{equation}
   \label{eq:uthm4_d1}
  - \frac{1}{2}\left(c_u - \frac{1}{3}\right)x^{3/2} + \left[1 + f_u(c_u)\right]x^2 \leq 0 ~~ \Rightarrow ~~ x \leq \left\{\frac{3c_u-1}{6\left[1+f_u(c_u)\right]}\right\}^2.
  \end{equation}

  Let's define $\d_3$ such that for $\d<\d_3$ \eqref{eq:uthm4_d1} holds; since $x$ is a decreasing function of $\d^{-1}$ for fixed $\g$ there exist a unique $\d_3$. We set $\d_0 = \min\left(\d_1,\d_2,\d_3\right)$ and conclude that if $c_u>1/3$, for fixed $\g>\g_0=4$ and $\e>0$ when $\d<\d_0$ as $\d\rightarrow 0$ \eqref{part1g_13} will remain negative and $2\Psi_{\max}\left(\widetilde{\lambda}^{\max}_{\g}+\epsilon,\delta,\rho\right)<0$.

  Having established a negative bound from above and the $\d_0$ for which it is valid, it remains to show that $~n \cdot 2\Psi_{\max}\left(\widetilde{\lambda}^{\max}_{\g}+\epsilon,\delta,\rho\right) \rightarrow -\infty$ as $(k,n,N) \rightarrow \infty$, which verifies an exponential decay to zero of the bound \eqref{ppmax} with $n$. This follows from the first term of the right hand side of \eqref{part1g_13}, giving a concluding bound $-n \cdot 4\left(c_u - 1/3\right)/\sqrt{\g^3}.$ For fixed $\g>\g_0$ and $\d<\d_0$ therefore
  \begin{multline*}
  Prob\left(U(k,n,N;A) > \widetilde{\mathcal{U}}^{\g}(\delta,\rho) + \e\right) \\ \le poly\left(n,\widetilde{\lambda}^{\max}_{\g}+\epsilon\right)\cdot\exp\left[-\frac{4n}{\sqrt{\g^3}}\left(c_u - \frac{1}{3}\right)\right].
  \end{multline*}
  The right hand side of which goes to zero as $n\rightarrow \infty$.\\

  \item[Part 2:]
  \begin{multline}
  \label{part2g_1}
  2\Psi_{\max}\left(\widetilde{\lambda}^{\max}_{\g}-\epsilon,\delta,\rho\right) = \left(1+\r\right)\log\left(\widetilde{\lambda}^{\max}_{\g}-\epsilon\right) - \r\log\left(\r\right)\\ + \r + 1 - \widetilde{\lambda}^{\max}_{\g} + \epsilon + \frac{2}{\delta}H\left(\delta\r\right),
  \end{multline}
  by substituting $\widetilde{\lambda}^{\max}_{\g}-\epsilon$ for $\lambda$ in the definition of $\Psi_{\max}\left(\lambda,\delta,\rho\right)$ in \eqref{eq:lmax}.

  Now letting $u = \widetilde{\lambda}^{\max}_{\g} - 1$ and substituting this in \eqref{part2g_1} and lower bounding the Shannon entropy term using the second bound of \eqref{eq:shannon_bounds} gives \eqref{part2g_2} below
  \begin{align}
  & 2\Psi_{\max}\left(\widetilde{\lambda}^{\max}_{\g}-\epsilon,\delta,\rho\right) \nonumber\\
  \label{part2g_2}
  & > \left(1+\r\right)\log\left(1+u-\epsilon\right) - \r\log\left(\r\right) + \r + 1 - (1+u) + \epsilon \qquad \nonumber\\
  & \quad + \frac{2}{\d}\left[-\d\r\log\left(\d\r\right) + \d\r - \d^2\r^2\right],\\
  \label{part2g_3}
  & = \log(1 + u - \e) + \r\log(1 + u-\e) - u + \epsilon + \r\log\left(\frac{1}{\d^2\r^3}\right) + 3\r \nonumber\\
  & \quad - 2\d\r^2,\\
  \label{part2g_4}
  & = \log(1 + u) + \log\left(1 - \frac{\e}{1 + u}\right) + \r\log(1 + u - \e) - u + \epsilon\nonumber\\
  & \quad + \r\log\left(\frac{1}{\d^2\r^3}\right) + 3\r - 2\d\r^2,\\
  \label{part2g_5}
  & > - u + u - \frac{1}{2}u^2 + \frac{1}{5}u^3 + \r\log\left(\frac{1}{\d^2\r^3}\right) + 3\r + \epsilon + \log\left(1 - \e\right)\nonumber\\
  & \quad + \r\log(1 + u - \e) - 2\d\r^2,\\
  \label{part2g_6}
  & = - \frac{1}{2}u^2 + \frac{1}{5}u^3 + \r\log\left(\frac{1}{\d^2\r^3}\right) + 3\r + \epsilon + \log\left(1 - \e\right) + \r\log(1 + u - \e)\qquad\nonumber\\
  & \quad - 2\d\r^2.
  \end{align}
  From \eqref{part2g_2} to \eqref{part2g_3} we expanded the $\left(1+\r\right)$ in the first term and simplified while from \eqref{part2g_3} to \eqref{part2g_4} we expanded the first logarithmic term. From \eqref{part2g_4} to \eqref{part2g_5} we bounded above $1/(1+u)$ using the bound of \eqref{eq:r_pfrac_series} and bounded below $\log(1 + u)$ using the following bound.
  \begin{equation}
   \label{eq:logp_lbound2}
  \log(1+x) \geq x - \frac{1}{2}x^2 + \frac{1}{5}x^3 \quad \forall x\in [0,0.92].
  \end{equation}
  From \eqref{part2g_5} to \eqref{part2g_6} we simplified.
  Let $x = 2\r\log\left(\frac{1}{\d^2\r^3}\right) + 6\r$ which means $u = \sqrt{x} + c_ux$. We simplify \eqref{part2g_6} and replace the second two terms by $x/2$ and $u$ in the first two terms by $\sqrt{x} + c_ux$ to get
  \begin{align}
  & 2\Psi_{\max}\left(\widetilde{\lambda}^{\max}_{\g}-\epsilon,\delta,\rho\right) \nonumber\\
  \label{part2g_7}
  & > - \frac{1}{2}\left(\sqrt{x} + c_ux\right)^2 + \frac{1}{5}\left(\sqrt{x} + c_ux\right)^3 + \frac{1}{2}x + \epsilon + \log\left(1 - \e\right)\qquad\nonumber\\
  & \quad + \r\log(1 + u - \e) - 2\d\r^2,\\
  \label{part2g_8}
  & = -\frac{1}{2}x - c_ux^{3/2} - \frac{1}{2}c_u^2x^2 + \frac{1}{5}x^{3/2} + \frac{3}{5}c_ux^2 + \frac{3}{5}c_u^2x^{5/2} + \frac{1}{5}c_u^3x^3 + \frac{1}{2}x\qquad\nonumber\\
  & \quad + \epsilon + \log\left(1 - \e\right) + \r\log(1 + u - \e) - 2\d\r^2,\\
  \label{part2g_9}
  & = \left(\frac{1}{5} - c_u\right)x^{3/2} + c_u\left(1 - \frac{1}{2}c_u\right)x^2 + \frac{3}{5}c_u^2x^{5/2} + \frac{1}{5}c_u^3x^3 + \r\log(1 + u - \e)\nonumber\\
  & \quad + \epsilon + \log\left(1 - \e\right) - 2\d\r^2.
  \end{align}

  From \eqref{part2g_7} to \eqref{part2g_8} we expanded the first two brackets and from \eqref{part2g_8} to \eqref{part2g_9} we simplified.  The dominant terms that does not go to zero as $\d \rightarrow 0$ are the terms with $x$ and their sum is positive for $c_u \leq 1/5$. Hence for fixed $\g$ there does not exist a $\d_0$ such that $2\Psi_{\max}\left(\widetilde{\lambda}^{\max}_{\g}-\epsilon,\delta,\rho\right) \leq 0.$ Thus
  \begin{multline*}
  Prob\left(U(k,n,N;A) > \widetilde{\mathcal{U}}^{\g}(\delta,\rho)-\e\right) \\ \le poly\left(n,\widetilde{\lambda}^{\max}_{\g}-\epsilon\right)\cdot\exp\left[2n\Psi_{\max}\left(\widetilde{\lambda}^{\max}_{\g}-\epsilon,\delta,\rho\right)\right],
  \end{multline*}
   and as $n \rightarrow \infty$ the right hand side of this does not go to zero.

\end{description}

Now \textbf{Part 1} and \textbf{Part 2} put together shows that $\widetilde{\mathcal{U}}^{\g}(\delta,\rho)$ is also a tight upper bound of $U(k,n,N;A)$ with overwhelming probability as the problem size grows in the regime prescribed for $\widetilde{\mathcal{U}}^{\g}(\delta,\rho)$ in Theorem \ref{thm:ULbctg}.

\end{proof}

\subsubsection{The lower bound, $\widetilde{\mathcal{L}}^{\g}(\delta,\rho_{\g}(\d))$}\label{sec:pLbctgthm}

\begin{proof}
Lets also define $$\widetilde{\lambda}^{\min}_{\g}(\delta,\rho) := 1 - \sqrt{2\rho\log\left(\frac{1}{\delta^2\rho^3}\right) + 6\r} + c_l\left[2\rho\log\left(\frac{1}{\delta^2\rho^3}\right) + 6\r\right].$$ This implies that $\widetilde{\mathcal{L}}^{\g}(\delta,\rho) = 1 - \widetilde{\lambda}^{\min}_{\g}(\delta,\rho)$ following from \eqref{eq:Lbctg}. Bounding $\widetilde{\mathcal{L}}^{\g}(\delta,\rho)$ above by $\widetilde{\mathcal{L}}^{\g}(\delta,\rho)+\e$ is equivalent to bounding $\widetilde{\lambda}^{\min}_{\g}$ below by $\widetilde{\lambda}^{\min}_{\g}-\epsilon$. We first establish that for a slightly looser bound, with $c_l > 1/3$, the exponent $\Psi_{\min}\left(\widetilde{\lambda}^{\min}_{\g}-\epsilon,\delta,\rho\right)$ is negative and then verify that when multiplied by $n$ it diverges to $-\infty$ as $n$ increases. We also show that for a slightly tighter bound, with $c_l<1/3$, the exponent $\Psi_{\min}\left(\widetilde{\lambda}^{\min}_{\g}+\epsilon,\delta,\rho\right)$ is bounded from below by zero, and hence the bound $\widetilde{\mathcal{L}}^{\g}(\delta,\rho)$ cannot be improved using the inequality \eqref{ppmin} from \cite{blanchard2011compressed}.  We show, in two parts that for $\g>\g_0$ fixed:
\begin{enumerate}
  \item $\exists ~\d_0, ~\epsilon > 0 ~\rm{and} ~c_l<1/3$ such that for $\d < \d_0, \Psi_{\min}\left(\widetilde{\lambda}^{\min}_{\g}-\epsilon,\delta,\rho\right) \leq 0;$
  \item $\nexists ~\d_0, ~\epsilon > 0 ~\rm{and} ~c_l\geq1/2$ such that for $\d < \d_0, \Psi_{\min}\left(\widetilde{\lambda}^{\min}_{\g}+\epsilon,\delta,\rho\right) \leq 0,$
\end{enumerate}
which are proven separately in the two parts as follows.

\begin{description}
  \item[Part 1:]
  \begin{multline}
  \label{part1g_1b}
  2\Psi_{\min}\left(\widetilde{\lambda}^{\min}_{\g}-\epsilon,\delta,\rho\right) = 2\mathrm{H}(\rho) + (1-\rho)\log\left(\widetilde{\lambda}^{\min}_{\g}-\epsilon\right)\\ + \rho\log(\rho) - \rho + 1 - \left(\widetilde{\lambda}^{\min}_{\g}-\epsilon\right)+ \frac{2}{\delta}\mathrm{H}(\delta\rho),	
  \end{multline}
  by substituting $\widetilde{\lambda}^{\min}_{\g}-\epsilon$ for $\lambda$ in \eqref{eq:lmin}. Let $l := 1 - \widetilde{\lambda}^{\min}_{\g}$ and bound the Shannon entropy functions from above using the first bound in \eqref{eq:shannon_bounds} which gives
  \begin{align}
  & 2\Psi_{\min}\left(\widetilde{\lambda}^{\min}_{\g}-\epsilon,\delta,\rho\right) \nonumber\\
  \label{part1g_2b}
  & < -2\r\log\left(\r\right) + 2\r + (1-\rho)\log\left[(1-l) - \e\right] + \rho\log\r - \rho + 1 - (1-l)  \nonumber\\
  & \quad + \e - 2\r\log\left(\d\r\right) + \frac{2}{\d}(\d\r) ,\\
  \label{part1g_3b}
  & = (1-\rho)\log\left(1-l-\e \right) + \rho\log\left(\frac{1}{\d^2\r^3}\right) + 3\rho + l + \e,\\
  \label{part1g_4b}
  & = l + \log(1-l) + \e + \log\left(1-\frac{\e}{1-l}\right) - \rho\log\left(1-l-\e\right)\qquad \nonumber\\
  & \quad + \rho\log\left(\frac{1}{\d^2\r^3}\right) + 3\rho, \\
  \label{part1g_5b}
  & < l + - l - \frac{1}{2}l^2 - \frac{1}{3}l^3 + \rho\log\left(\frac{1}{\d^2\r^3}\right) + 3\rho - \rho\log\left(1-l-\e\right) \nonumber\\
  & \quad + \log(1-\e) + \e, \\
  \label{part1g_6b}
  & < - \frac{1}{2}l^2 - \frac{1}{3}l^3 + \rho\log\left(\frac{1}{\d^2\r^3}\right) + 3\rho - \rho\log\left(1-l-\e\right) - \e + \e.
  \end{align}
  We simplified from \eqref{part1g_2b} to \eqref{part1g_3b} and from \eqref{part1g_3b} to \eqref{part1g_4b} we expanded the first logarithmic term. From \eqref{part1g_4b} to \eqref{part1g_5b} we bounded $1/(1-l)$ below and $\log(1-l)$ above using \eqref{eq:r_mfrac_series} and the third bound of \eqref{eq:logm_ubound} respectively. From \eqref{part1g_5b} to \eqref{part1g_6b} we simplified and bounded above $\log(1-\e)$ using the first bound of \eqref{eq:logm_ubound}.

  Let $x = 2\r\log\left(\frac{1}{\d^2\r^3}\right) + 6\r$ which means $l = \sqrt{x} - c_lx$. We simplify \eqref{part1g_6b} and replace the second two terms by $\frac{1}{2}x$ and $l$ in the first two terms by $\sqrt{x} - c_lx$ to get
  \begin{align}
  & 2\Psi_{\min}\left(\widetilde{\lambda}^{\min}_{\g}-\epsilon,\delta,\rho\right)\nonumber\\
  \label{part1g_7b}
  & < - \frac{1}{2}\left(\sqrt{x} - c_lx\right)^2 - \frac{1}{3}\left(\sqrt{x} - c_lx\right)^3 + \frac{1}{2}x - \rho\log\left(1-l-\e\right)\\
  \label{part1g_8b}
  & = -\frac{1}{2}x + c_lx^{3/2} - \frac{1}{2}c_l^2x^2 - \frac{1}{3}x^{3/2} + c_lx^2 - c_l^2x^{5/2} + \frac{1}{3}c_l^3x^3 + \frac{1}{2}x \nonumber\\
  & \quad - \rho\log\left(1-l-\e\right),\\
  \label{part1g_9b}
  & = - \left(\frac{1}{3} - c_l\right)x^{3/2} + c_lx^2 - \frac{1}{2}c_l^2x^2 - c_l^2x^{5/2} + \frac{1}{3}c_l^3x^3\nonumber\\
  & \quad - \rho\log\left(1-l-\e\right).
  \end{align}
  From \eqref{part1g_7b} to \eqref{part1g_8b} we expanded the first two brackets and from \eqref{part1g_8b} to \eqref{part1g_9b} we simplified. Substituting $1/\left[\g\log\left(1/\d\right)\right]$ for $\r$ in the expression for $x$ we have $x = 4/\g + g(\r)$ where $g(\r) = 6\r\log\left(1/\r\right) + 6\r$ and goes to zero with $\d$. We make the same argument as in \textbf{Part 1} of the proof for $\widetilde{\mathcal{U}}^{\g}(\delta,\rho_{\g}(\d))$ in Section \ref{sec:pLbctgthm}, that is for $\g > 4$ we can define $\d_1$ such that for $\d<\d_1, ~x<1$ and we can upper bound $x^3$ by $x^{2}$ since $x^{2} > x^{2+j}$ for $j>0$ when $x<1$. The last term in \eqref{part1g_9b} goes to zero with $\d$, so we can define $\d_2$ such that for $\d<\d_2$ we can bound this term above by $x^2$ which is a constant. We split the first term of \eqref{part1g_9b} into half and drop the two $c_l^2$ terms because they are negative. Let $f_l(c_l) = c_l + c_l^3/3$, which is positive for all $c_l>0$, using the above we upper bound \eqref{part1g_9b} as follows.
  \begin{align}
  & 2\Psi_{\min}\left(\widetilde{\lambda}^{\min}_{\g}-\epsilon,\delta,\rho\right)\nonumber\\
  \label{part1g_10b}
  & < - \frac{1}{2}\left(\frac{1}{3} - c_l\right)x^{3/2}  - \frac{1}{2}\left(\frac{1}{3} - c_l\right)x^{3/2} + f_l(c_l)x^2 + x^2,\\
  \label{part1g_11b}
  & < -\frac{4}{\sqrt{\g^3}}\left(\frac{1}{3} - c_l\right) - \frac{1}{2}\left(\frac{1}{3} - c_l\right)x^{3/2} + \left[1+f_l(c_l)\right]x^2.
  \end{align}
  From \eqref{part1g_10b} to \eqref{part1g_11b} we use the fact that $-x^{3/2} < -8/\sqrt{\g^3}$ as shown in Section \ref{sec:pLbctgthm}. For \eqref{part1g_11b} to be negative all we need is for $c_l<1/3$ and the sum of the last two terms to be non positive, that is:
  \begin{equation}
   \label{eq:lthm4_d1}
    - \frac{1}{2}\left(\frac{1}{3} - c_l\right)x^{3/2} + \left[1+f_l(c_l)\right]x^2 \leq 0 \quad \Rightarrow \quad x \leq \left\{\frac{1-3c_l}{6\left[1+f_l(c_l)\right]}\right\}^2.
  \end{equation}

  Let's define $\d_3$ such that for $\d<\d_3$ \eqref{eq:lthm4_d1} holds; since $x$ is a decreasing function of $\d^{-1}$ for fixed $\g$ there exist a unique $\d_3$. We set $\d_0 = \min\left(\d_1,\d_2,\d_3\right)$ and conclude that if $c_l<1/3$, for fixed $\g>\g_0=4$ and $\e>0$ when $\d<\d_0$ as $\d\rightarrow 0$ \eqref{part1g_11b} will remain negative and $2\Psi_{\min}\left(\widetilde{\lambda}^{\min}_{\g}-\epsilon,\delta,\rho\right)<0$.

  Having established a negative bound from above and the $\d_0$ for which it is valid, it remains to show that $~n \cdot 2\Psi_{\min}\left(\widetilde{\lambda}^{\min}_{\g}-\epsilon,\delta,\rho\right) \rightarrow -\infty$ as $(k,n,N) \rightarrow \infty$, which verifies an exponential decay to zero of the bound \eqref{ppmin} with $n$. This follows from the first term of the right hand side of \eqref{part1g_11b} giving a concluding bound $-n\cdot 4 \left(1/3 - c_l\right)/\sqrt{\g^3}.$ For $\g>\g_0$ and $\d<\d_0$ therefore
  \begin{multline*}
  Prob\left(L(k,n,N;A) > \widetilde{\mathcal{L}}^{\g}(\delta,\rho) + \e\right) \\ \le poly\left(n,\widetilde{\lambda}^{\min}_{\g}+\epsilon\right)\cdot\exp\left[-\frac{4n}{\sqrt{\g^3}}\left(\frac{1}{3} - c_l\right)\right].
  \end{multline*}
  The right hand side of which goes to zero as $n\rightarrow \infty$.\\

  \item[Part 2:]
  \begin{multline}
  \label{part2g_1b}
  2\Psi_{\min}\left(\widetilde{\lambda}^{\min}_{\g}+\epsilon,\delta,\rho\right) = 2\mathrm{H}(\rho) + (1-\rho)\log\left(\widetilde{\lambda}^{\min}_{\g}+\epsilon\right)\\ + \rho\log(\rho) - \rho + 1 - \left(\widetilde{\lambda}^{\min}_{\g}+\epsilon\right)+ \frac{2}{\delta}\mathrm{H}(\delta\rho),	
  \end{multline}
  by substituting $\widetilde{\lambda}^{\min}_{\g}+\epsilon$ for $\lambda$ in \eqref{eq:lmin}. Let $l := 1 - \widetilde{\lambda}^{\min}_{\g}$ and bound the Shannon entropy function from below using the second bound in \eqref{eq:shannon_bounds} to give
  \begin{align}
  & 2\Psi_{\min}\left(\widetilde{\lambda}^{\min}_{\g}+\epsilon,\delta,\rho\right) \nonumber\\
  \label{part2g_2b}
  & > 2 \left[-\r\log\r + \r - \r^2\right] + (1-\rho)\log\left[(1-l) + \e\right] + \rho\log\r - \rho  \nonumber\\
  & \quad + 1 - \left(1-l\right) - \e + \frac{2}{\d}\left[-\r\log\left(\d\r\right) + \d\r - \d^2\r^2\right],\\
  \label{part2g_3b}
  & = -2\r\log\r + 2\r - 2\r^2 + \log\left(1 - l + \e \right) - \rho\log\left(1 - l + \e \right) + \rho\log\r - \rho \nonumber\\
  & \quad + 1 - 1 + l - \e - 2\r\log\left(\d\r\right) + 2\r - 2\d\r^2,\\
  \label{part2g_4b}
  & = \log\left(1-l\right) + \log\left(1 + \frac{\e}{1-l}\right) + l - \e - \rho\log\left(1 - l + \e \right) + \rho\log\left(\frac{1}{\d^2\r^3}\right)\nonumber\\
  & \quad + 3\rho - 2(1+\d)\r^2,\\
  \label{part2g_5b}
  & > -l - \frac{1}{2}l^2 - \frac{1}{2}l^3 + l + \rho\log\left(\frac{1}{\d^2\r^3}\right) + 3\rho + \log\left(1 + \e\right) - \e\nonumber\\
  & \quad - \rho\log\left(1-l+\e\right) - 2(1-\d)\r^2, \\
  \label{part2g_6b}
  & > - \frac{1}{2}l^2 - \frac{1}{2}l^3 + \rho\log\left(\frac{1}{\d^2\r^3}\right) + 3\rho + \e - \frac{1}{2}\e^2 - \e - \rho\log\left(1-l+\e\right)\nonumber\\
  & \quad - 2(1-\d)\r^2.
  \end{align}
  From \eqref{part2g_2b} to \eqref{part2g_3b} we expanded brackets and simplified. From \eqref{part2g_3b} to \eqref{part2g_4b} we expanded $\log\left(1 - l + \e \right)$ and simplified. From \eqref{part2g_4b} to \eqref{part2g_5b} we bounded from below $1/(1-l)$ using \eqref{eq:r_mfrac_series} and using the bound of \eqref{eq:logm_lbound} we also bounded from below $\log(1 - l)$. Then from \eqref{part2g_5b} to \eqref{part2g_6b} we simplified and bounded from below $\log(1+\e)$ using \eqref{eq:logp_lbound}.

  Let $x = 2\r\log\left(\frac{1}{\d^2\r^3}\right) + 6\r$ which means $l = \sqrt{x} - c_lx$. We simplify \eqref{part2g_6b} and replace the second two terms by $x/2$ and $l$ in the first two terms by $\sqrt{x} - c_lx$ to get
  \begin{align}
  & 2\Psi_{\min}\left(\widetilde{\lambda}^{\min}_{\g}+\epsilon,\delta,\rho\right)\nonumber\\
  \label{part2g_7b}
  & > - \frac{1}{2}\left(\sqrt{x} - c_lx\right)^2 - \frac{1}{2}\left(\sqrt{x} - c_lx\right)^3 + \frac{1}{2}x - \rho\log\left(1-l+\e\right)\nonumber\\
  & \quad - 2(1-\d)\r^2 - \frac{1}{2}\e^2,\\
  \label{part2g_8b}
  & = -\frac{1}{2}x + c_lx^{3/2} - \frac{1}{2}c_l^2x^2 - \frac{1}{2}x^{3/2} + \frac{3}{2}c_lx^2 - \frac{3}{2}c_l^2x^{5/2} + \frac{1}{2}c_l^3x^3 + \frac{1}{2}x \qquad\nonumber\\
  & \quad - \rho\log\left(1-l+\e\right) - 2(1-\d)\r^2 - \frac{1}{2}\e^2,\\
  \label{part2g_9b}
  & = \left(c_l - \frac{1}{2}\right)x^{3/2} + \frac{1}{2}c_l\left(3 - c_l\right)x^2 - \frac{3}{2}c_l^2x^{5/2} + \frac{1}{2}c_l^3x^3 - \rho\log\left(1-l+\e\right) \nonumber\\
  & \quad - 2(1-\d)\r^2 - \frac{1}{2}\e^2.
  \end{align}
  From \eqref{part2g_7b} to \eqref{part2g_8b} we expanded the first two brackets and simplified from \eqref{part2g_8b} to \eqref{part2g_9b}. The dominant terms that does not go to zero as $\d \rightarrow 0$ are the terms with $x$ and their sum is positive if $c_l \geq 1/2$ and $x<1$. We established in the earlier parts of this proof of Theorem \ref{thm:ULbctg} that if $\g>4$ we will have $x<1$ as $\d\rightarrow 0$. Hence we conclude that for fixed $\g>\g_0=4$ and $\e>0$ there does not exist a $\d_0$ such that \eqref{part2g_9b} is negative and $2\Psi_{\min}\left(\widetilde{\lambda}^{\min}_{\g}+\epsilon,\delta,\rho\right) \leq 0$ as $\d\rightarrow 0$.
  Thus
  \begin{multline*}
  Prob\left(L(k,n,N;A) > \widetilde{\mathcal{L}}^{\g}(\delta,\rho)-\e\right) \\ \le poly\left(n,\widetilde{\lambda}^{\min}_{\g}+\epsilon\right)\cdot\exp\left[2n\Psi_{\min}\left(\widetilde{\lambda}^{\min}_{\g}+\epsilon,\delta,\rho\right)\right],
  \end{multline*}
  and as $n \rightarrow \infty$ the right hand side of this does not go to zero.

\end{description}
Now \textbf{Part 1} and \textbf{Part 2} put together shows that $\widetilde{\mathcal{L}}^{\g}(\delta,\rho)$ is also a tight bound of $L(k,n,N;A)$ with overwhelming probability as the sample size grows in the regime prescribed for $\widetilde{\mathcal{L}}^{\g}(\delta,\rho)$ in Theorem \ref{thm:ULbctg}.

\end{proof}

\subsection{Corollary \ref{cor_ULrgd}}\label{sec:pULrgdcor}

\begin{proof}

 We prove Corollary \ref{cor_ULrgd} in two parts, first proving the case for $\widetilde{\mathcal{U}}^{\g}(\delta,\rho_{\g}(\d))$ and then that of $\widetilde{\mathcal{L}}^{\g}(\delta,\rho_{\g}(\d)).$

\begin{description}
  \item[Part 1:]
  From \eqref{eq:Ubctg}, for $\r = \r_{\g}(\d) = \left(\g \log\left(\frac{1}{\d}\right)\right)^{-1},$ we have
  \begin{align}
  \label{eq:urgd1}
  \widetilde{\mathcal{U}}^{\g}(\delta,\rho_{\g}(\d)) &= \sqrt{2\rho\log\left(\frac{1}{\delta^2\rho^3}\right) + 6\r} + c_u\left[2\rho\log\left(\frac{1}{\delta^2\rho^3}\right) + 6\r\right],\\
  \label{eq:urgd2}
  & = \sqrt{2\rho\log\left(\frac{1}{\delta^2\rho^3}\right) + 6\r} + 2c_u\rho\log\left(\frac{1}{\delta^2\rho^3}\right) + 6c_u\r\\
  \label{eq:urgd3}
  & = \sqrt{4\rho\log\left(\frac{1}{\delta}\right) + 6\rho\log\left(\frac{1}{\rho}\right) + 6\r} + 4c_u\rho\log\left(\frac{1}{\delta}\right)\nonumber \\
  & \quad + 6c_u\rho\log\left(\frac{1}{\rho}\right) + 6c_u\r,\\
  \label{eq:urgd4}
  & = \sqrt{\frac{4}{\g} + 6\rho\log\left(\frac{1}{\rho}\right) + 6\r} + \frac{4c_u}{\g} + 6c_u\rho\log\left(\frac{1}{\rho}\right) \nonumber\\
  & \quad + 6c_u\r.
  \end{align}

  From \eqref{eq:urgd1} to \eqref{eq:urgd2} we expanded the square brackets while from \eqref{eq:urgd2} to \eqref{eq:urgd3} we separated the terms explicitly involving $\d$ from the rest. From \eqref{eq:urgd3} to \eqref{eq:urgd4} we substituted $1/\left[\g \log\left(1/\d\right)\right]$ for $\r$ in the terms explicitly involving $\d$ and simplified.

  Now using the fact that $\lim_{\delta\rightarrow 0} \r\log\left(1/\r\right) = 0$ and $\lim_{\delta\rightarrow 0} \r = 0$ we have
  \begin{equation*}
  \lim_{\delta\rightarrow 0} \widetilde{\mathcal{U}}^{\g}(\delta,\rho_{\g}(\d)) = \frac{2}{\sqrt{\g}} + \frac{4c_u}{\g},
  \end{equation*}
 hence concluding the proof for $\widetilde{\mathcal{U}}^{\g}(\delta,\rho_{\g}(\d)).$\\

  \item[Part 2:]
  From \eqref{eq:Lbctg}, for $\r = \r_{\g}(\d) = \left(\g \log\left(\frac{1}{\d}\right)\right)^{-1},$ we have
  \begin{align}
  \label{eq:lrgd1}
  \widetilde{\mathcal{L}}^{\g}(\delta,\rho_{\g}(\d)) &= \sqrt{2\rho\log\left(\frac{1}{\delta^2\rho^3}\right) + 6\r} - c_l\left[2\rho\log\left(\frac{1}{\delta^2\rho^3}\right) + 6\r\right],\\
  \label{eq:lrgd2}
  & = \sqrt{2\rho\log\left(\frac{1}{\delta^2\rho^3}\right) + 6\r} - 2c_l\rho\log\left(\frac{1}{\delta^2\rho^3}\right) - 6c_l\r,\\
  \label{eq:lrgd3}
  & = \sqrt{4\rho\log\left(\frac{1}{\delta}\right) + 6\rho\log\left(\frac{1}{\rho}\right) + 6\r} - 4c_l\rho\log\left(\frac{1}{\delta}\right)\nonumber \\
  & \quad - 6c_l\rho\log\left(\frac{1}{\rho}\right) - 6c_l\r,\\
  \label{eq:lrgd4}
  & = \sqrt{\frac{4}{\g} + 6\rho\log\left(\frac{1}{\rho}\right) + 6\r} - \frac{4c_l}{\g} - 6c_l\rho\log\left(\frac{1}{\rho}\right)\nonumber\\
  & \quad - 6c_l\r.
  \end{align}

  From \eqref{eq:lrgd1} to \eqref{eq:lrgd2} we expanded the square brackets while from \eqref{eq:lrgd2} to \eqref{eq:lrgd3} we separated the terms explicitly involving $\d$ from the rest. Then from \eqref{eq:lrgd3} to \eqref{eq:lrgd4} we substituted $1/\left[\g \log\left(1/\d\right)\right]$ for $\r$ in the terms explicitly involving $\d$ and simplified.

  Now using the fact that $\lim_{\delta\rightarrow 0} \r\log\left(1/\r\right) = 0$ and $\lim_{\delta\rightarrow 0} \r = 0$ we have
  \begin{equation*}
  \lim_{\delta\rightarrow 0} \widetilde{\mathcal{L}}^{\g}(\delta,\rho_{\g}(\d)) = \frac{2}{\sqrt{\g}} - \frac{4c_l}{\g},
  \end{equation*}
 hence concluding the proof for $\widetilde{\mathcal{U}}^{\g}(\delta,\rho_{\g}(\d)).$
\end{description}

\textbf{Part 1} and \textbf{Part 2} combined concludes the proof for Corollary \ref{cor_ULrgd}.
\end{proof}

%
%\bibliographystyle{elsarticle-num}
%\bibliography{asymptotic_rics}

\end{document}